\documentclass[a4paper,12pt]{article}
\usepackage{graphicx}
\usepackage{amsmath,amsthm}
\usepackage{amsfonts,amssymb}
\usepackage{amscd}
\usepackage{color}
\usepackage[cp1251]{inputenc}
\usepackage[matrix,arrow,curve]{xy}

\newtheorem{lemma}{Lemma}
\newtheorem{prop}[lemma]{Proposition}
\newtheorem{theorem}[lemma]{Theorem}
\newtheorem{coro}[lemma]{Corollary}
\newtheorem{rema}[lemma]{{\rm\bf Remark}}

\newtheorem{thm}{Theorem}
\newtheorem{Def}[thm]{Definition}

\numberwithin{lemma}{section}

\newtheorem{ex}{Example}

\newcommand{\Tot}{{\rm Tot}}
\newcommand{\Hom}{{\rm Hom}}
\newcommand{\Homk}{{\rm Hom}_{\kk}}
\newcommand{\charr}{{\rm char\,}}
\newcommand{\kk}{{\bf k}}
\newcommand{\mm}{{\bf m}}
\newcommand{\hh}{{\bf h}}
\newcommand{\bh}{{\bf \bar{h}}}

\DeclareMathOperator{\Ker}{{Ker}}

\DeclareMathOperator{\Aut}{{Aut}}
\renewcommand{\Bar}{{\rm Bar}}
\newcommand{\Id}{{\rm Id}}
\newcommand{\Ext}{{\rm Ext}}
\newcommand{\Mod}{{\rm Mod}}
\newcommand{\Der}{{\rm Der}}

\DeclareMathOperator{\HH}{{HH}}
\DeclareMathOperator{\Ho}{{H}}
\renewcommand{\Im}{{\rm Im\,}}

\newcommand{\NN}{{\mathcal N}}
\newcommand{\XX}{{\mathcal X}}
\newcommand{\RR}{{\mathcal R}}
\newcommand{\AAA}{{\mathcal A}}

\newcommand{\ee}{\epsilon}
\newcommand{\ve}{\varepsilon}

\renewcommand{\le}{\leqslant}
\renewcommand{\ge}{\geqslant}
\newcommand{\upper}[1]{\left\lceil#1\right\rceil}

\newenvironment{Ex}{\begin{ex}
\rm }{\flushright
$\Box$\end{ex}}

\newenvironment{Proof}[1][Proof.]{\begin{trivlist}
\item[\hskip \labelsep {\bfseries #1}]}{\flushright
$\Box$\end{trivlist}}

\numberwithin{equation}{section}
\begin{document}
\title{Hochschild cohomology of a smash product in the nonsemisimple case}
\author{Eduardo do Nascimento Marcos \footnote{The first named author  supported by  CNPQ (grant number  300858/2004-3), and also by a tematic project of Fapesp (Project number:  
2014/09310-5    )}\ \ and Yury Volkov
\footnote{The second named author was supported by a pos-doc scholarship of Fapesp process number 2014/19521-3}}
\date{}
\maketitle
\begin{abstract}
In this paper we explore the relations between the Hochschild cohomology of an algebra over some field and the Hochschild cohomology of its smash product with a finite group. 
Basically we are concentrated on the case where the group under consideration is an extension of a cyclic $p$-group by some $p'$-group, where $p$ is the characteristic of the ground field.
\end{abstract}

\section{Introduction}

Hochschild cohomology is a very important derived invariant of an algebra. It contains many important information about its derived category. It is also related to various 
properties, for instance  it is connected with its simple connectedness, the fundamental group of a presentation etc. In the commutative case it tells 
a lot about the
smoothness. It is also important to know how the Hochschild cohomology of related algebra are related. This work looks at some aspects of such 
relations, namely, at the relation between 
the Hochschild cohomology of an algebra with a group action and of its skew group algebra.

Let $A$ be a finite-dimensional unital algebra over a field $\kk$ and $G$ be a finite group acting on $A$. Then we can define a smash product of $A$ 
and $\kk G$ (sometimes it is called a skew group ring), which is denoted by $A\#\kk G$.
Let $M$ be an $A\#\kk G$-algebra. It is well known that in this case $G$ acts on $\HH^*(A,M)$. Moreover, there is a spectral sequence of algebras
$$
E_2^{i,j}=\Ho^i(G,\HH^j(A,M))\Rightarrow \HH^{i+j}(A\#\kk G,M).
$$
If $\charr\kk$ does not divide $|G|$, then this sequence gives an isomorphism of algebras $\HH^*(A,M)^G\cong\HH^*(A\#\kk G,M)$. All these facts were 
proved in many works, see, for example, \cite{Sanada}, \cite{Stefan},  and \cite{CibRed}.

The facts mentioned above give a good relation between the Hochschild cohomology of $A$ and the Hochschild cohomology of $A\#\kk G$ in the case where the group algebra $\kk G$ is semisimple.
On the other hand, there is no good connection in the $p$-modular case. Among works in this direction the work \cite{Wit} have to be mentioned. Its authors construct a bimodule projective resolution of $A\#\kk G$ using bimodule projective resolutions of $A$ and $\kk G$ satisfying some properties. The authors of \cite{CibRed} do the same using a bimodule projective resolution of $A$ and a $\kk G$-projective resolution of a trivial module. They use this construction to obtain a spectral sequence connecting the Hochschild homology of $A$ and $A\#\kk G$. In this paper we use a similar construction to explore the Hochschild cohomology.

An attempt to obtain some results in the $p$-modular case was made by the second author in \cite{Volk5}. The results obtained there do not give a good general answer, but they allow to connect the Hochschild cohomology of $A$ and the Hochschild cohomology of $A\#\kk C_3$ for some specific algebra $A$ over a field of characteristic $3$.

This work was motivated by the following situation. Let $\RR$ be an algebra, $T\RR$ be its trivial extention by $D\RR$ and $\RR_n$ be an 
algebra $\hat{\RR}/\nu_{\infty}^n$, where $\hat{\RR}$ is a repetetive algebra of $\RR$ and $\nu_{\infty}$ its Nakayama automorphism. Then it is well 
known that $T\RR$ is Morita equivalent to $\RR_n\#\kk C_n$, where $C_n$ is generated by a Nakayama automorphism of $\RR_n$. To see 
this one can apply the duality theorem of Cohen-Montgomery (see \cite{CohMon}) to the fact that $\RR_n$ is a smash product of $T\RR$ with the 
cyclic group $C_n$. Using results of \cite{Volk6} we see that
any Nakayama automorphism acts trivially on Hochschild cohomology with coefficients in the regular module and so $\HH^*(\RR_n)$ is a subalgebra of 
$\HH^*(T\RR)$ if $\charr\kk\nmid n$. So the description of $\HH^*(\RR_n)$ can be easily obtained from  the description of $\HH^*(T\RR)$. 
The examples which we know indicates that  almost same fact has to be true in the case $\charr\kk\mid n$. This motivates us to study this case and 
obtain a desirable result (see Theorem \ref{TR}).

We recall that an algebra $B= A/I$ is called a singular extension of $A$ by $I$ if $I^2=0$. The current work is mainly devoted to the case where $G$ is a trivial extension of a cyclic $p$-group by some $p'$-group, where $p=\charr\kk$. 
Moreover, the most powerful results were obtained in the case where the spectral sequence mentioned above is $(3,2)$-degenerated, i.e. its third page has only two nonzero columns.
In the most general terms, our result in this case is as follows. There is a subalgebra $\AAA_M^G$ and an ideal $W_M^G\subset\AAA_M^G$ of $\HH^*(A,M)^G$ such that
$$
\HH^*(A,M)^G/W_M^G\cong (\AAA_M^G/W_M^G)[x,D]/\langle x^2-a\rangle,
$$
where $D$ is some graded derivation and $a$ is some element of $\AAA_M^G/W_M^G$.  
 It is also true that $\HH^*(AG,M)$ is a singular extension of the algebra $\AAA_M^G$ by some $\AAA_M^G$-bimodule $X_M^G$. And moreover, there are filtrations
$$0=X_{q-1}\subset\cdots\subset X_0=X_M^G\mbox{ and }0=Y_{q-1}\subset\cdots\subset Y_0=(\AAA_M^G/W_M^G)[-1]$$ such that $Y_{i}/Y_{i+1}\cong X_{q-2-i}/X_{q-1-i}$ for $0\le i\le q-2$. These results are contained in Theorem \ref{Gstabdes} and Theorem \ref{filtering}. Note also that $D=0$ if $M=A$ or $M=A\#\kk G$, and, if additionally $p\not=2$, then $a=0$ for the same $M$.

Among other results of this work we want to mention Theorem \ref{degeneration} and Theorem \ref{dims}. Theorem \ref{degeneration} gives a way to check 
if  the spectral sequence considered in this work degenerates (i.e. some of its pages have only finite number of nonzero columns). During 
its proof we obtain also the following interesting fact. If $G$ is a finite group, then there is some integer $S$ such that for any $\kk G$-module $L$ 
we have $\Ho^*(G,L)=\Ho^*(G,\kk)\smile\oplus_{i=0}^S\Ho^i(G,L)$. Theorem \ref{dims} establishes some relations between $\dim_{\kk}\HH^n(A,M)^G$ and 
$\dim_{\kk}\HH^n(AG,M)$ in the case where $G$ is an extension of cyclic $p$-group by some $p'$-group. It says that $\dim_{\kk}\HH^n(A,M)^G\le 
\dim_{\kk}\HH^n(AG,M)$ for any $n$ and $\dim_{\kk}\HH^1(A,M)^G= \dim_{\kk}\HH^1(AG,M)$ iff $\dim_{\kk}\HH^n(A,M)^G= \dim_{\kk}\HH^n(AG,M)$ for any 
$n\ge 0$.

\section{Hochschild cohomology}\label{HH}

In this section we recall the definition of Hochschild cohomology and related notions. All algebras in this paper are unital over some fixed field $\kk$. We write simply $\otimes$ instead of $\otimes_{\kk}$.
All modules in this paper, if otherwise is not stated, are left modules and we will denote by $\Mod A$ the category of modules over the algebra  
$A$.

\begin{Def}[{\bf Complexes}] 
An $A$-complex is a $\mathbb{Z}$-graded $A$-module $P$ with a differential of degree $-1$, i.e. an $A$-module $P$ 
with some 
fixed 
$A$-module direct sum decomposition $P=\oplus_{n\in \mathbb{Z}}P_n$ and an $A$-module homomorphism $d_P:P\rightarrow P$ such that $d_P(P_n)\subset 
P_{n-1}$ and $d_P^2=0$.  Let $d_{P,n}$ denote $d_P|_{P_{n}}$. The $n$-th homology of $P$ is the vector space $\Ho_n(P)=(\Ker d_{P,n})/(\Im d_{P,n+1})$. 
An $A$-complex $P$ is called acyclic if $\Ho_n(P)=0$ for all $n\in\mathbb{Z}$. A morphism of $A$-complexes is a homomorphism of $A$-modules that 
respects the respective grading and differential. A complex is called positive if $P_n=0$ for $n<0$.

If $P$ is a complex, then we denote by $P[t]$ the complex, which equals $P$ as an $A$-module, with grading $P[t]_n=P_{t+n}$ and 
differential defined as $d_{P[t]}= (-1)^t d_{P} $.
\end{Def}

Given an $A$-complex $P$ and $M\in\Mod A$, $\Hom_A(P,M)$ denotes the $\kk$-complex $\oplus_{n\le 0}\Hom(P_{-n},M)$ with differential $d_{\Hom_A(P,M),n}=\Hom_A(d_{P,-1-n},M)$.
Given $M\in\Mod A$, an $A$-resolution of $M$ is a pair $(P,\mu_P)$, where $P$ is a positive $A$-complex.
and $\mu_P:P_0\rightarrow M$ is a homomorphism of $A$-modules, which induces an isomorphism $\Ho_0(P)\cong M$.

From now on if $A$ is an algebra we denote by $\mu_A: A\otimes A \to A$ the multiplication map.

Given an algebra $A$, we denote by $A^e$ the algebra $A\otimes A^{\rm op}$, which is called the enveloping algebra of $A$. Let $\Bar(A)$ be the
positive $A^e$-complex with $n$-th member $\Bar_n(A)=A^{\otimes(n+2)}$ for $n\ge 0$ and the differential $d_{\Bar(A)}$ defined by the equality
 $$
d_{\Bar(A)}(a_0\otimes\cdots\otimes a_{n+1})=\sum\limits_{i=0}^n(-1)^ia_0\otimes\cdots\otimes a_{i-1}\otimes a_ia_{i+1}\otimes 
a_{i+2}\otimes\cdots\otimes a_{n+1}
 $$
for $n>0$ and $a_i\in A$ ($0\le i\le n+1$). 
Then $(\Bar(A),\mu_A)$, is a projective 
$A^e$-resolution of $A$, which is called bar resolution.

\begin{Def}[{\bf Hochschild Cohomology}]
The Hochschild cohomology of the algebra $A$ with coefficients in the $A^e$-module $M$ is the homology of the complex $C(A,M)=\Hom_{A^e}(\Bar(A),M)$. 
\end{Def}

We write $C^n(A,M)$ instead of $C_{-n}(A,M)$ and $\delta_M^n$ instead of $d_{C(A,M),-1-n}$.
Note that $C^0(A,M)\simeq M$ and $C^n(A,M)\simeq\Homk(A^{\otimes n},M)$.\\
Given $f\in C^n(A,M)$, we introduce the notation
$$
 \delta_n^i(f)(a_1\otimes\cdots\otimes a_{n+1}):=\begin{cases}
 a_1f(a_2\otimes\cdots\otimes a_{n+1}),&\mbox{if $i=0$},\\
  (-1)^if(a_1\otimes \cdots \otimes a_ia_{i+1} \otimes \cdots\otimes
a_{n+1}),&\mbox{if $1\le i\le n$},\\
(-1)^{n+1}f(a_1\otimes \cdots \otimes a_n) a_{n+1},&\mbox{if $i=n+1$}.
 \end{cases}
 $$
Then $\delta_M^n=\sum\limits_{i=0}^{n+1}\delta_n^i$. We have $\HH^n(A,M)=(\Ker \delta_M^n)/(\Im \delta_M^{n-1})$. From now on we write $C^n(A)$ and 
$\HH^n(A)$ instead of $C^n(A, A)$ and $\HH^n(A, A)$ respectively. It follows from the definition that  $\HH^n(A, M)= \Ext^n_{A^e}(A,M).$

\subsection{Tensor Product of Complexes, Cup Products}

We continue by defining tensor product of complexes and various ways cup products, which are  graded algebra products on the cohomology.

Let $\HH^*(A,M)$ denote the direct sum $\oplus_{n\ge 0}\HH^n(A,M)$.
\begin{Def}[ {\bf $A$-algebra}]
If $M$ is an $A$-bimodule and there is an $A$-bimodule map,\\ $M\otimes_A M \to M$, which defines a $\kk$-algebra structure on $M$, then $M$  will be 
called an $A$-algebra.
\end{Def}

\begin{Def}[{\bf Cup Product}]

Assume that $M$ is an $A$-algebra, $f\in C^n(A,M)$ and  \\ $g\in C^m(A,M)$.
 Then $f\smile g \in 
C^{n+m}(A,M)=\mathrm{Hom}_{\kk}(A^{\otimes (n+m)}, M)$ is defined by the formula
$$
(f\smile g)(a_1\otimes \cdots\otimes
a_{n+m}):=f(a_1\otimes\cdots\otimes a_n)\cdot
g(a_{n+1}\otimes\cdots\otimes a_{n+m}).
$$
If $M$ has a unity, then $1_M\in C^0(A,M)$ is a unity for the multiplication $\smile$.

In this case  $\smile$ induces a graded $\kk$-algebra structure on $\HH^*(A,M)$
$$
\smile: \HH^n(A,M) \otimes \HH^m(A,M) \longrightarrow \HH^{n+m}(A,M). 
$$
This product is know as cup product.
\end{Def}

\begin{Def}[\bf Tensor Product of Complexes]
Given $M$  an $A^{op}$-complex and $N$ an $A$-complex, we define the tensor product complex $M\otimes_A N$ by the equality $(M\otimes_A
N)_n=\sum_{i+j=n}M_i\otimes_A N_j$.  The differential $d_{M\otimes N}$ is defined by the equality $d_{M\otimes_A N}(u\otimes v)=d_M(u)\otimes_A
v+(-1)^iu\otimes_A d_N(v)$ for $u\in M_i$, $v\in N$. 
\end{Def}

\begin{Def}[\bf Comultiplication on a Projective Resolution]
Given a projective $A^e$-resolution $(P,\mu_P)$ of $A$, a morphism  $\Delta_P:P\rightarrow P\otimes_A P$ of 
$A^e$-complexes will be called a comultiplication on $P$ if $\mu_A(\mu_P\otimes\mu_P)\Delta_P=\mu_P$.
\end{Def}

If $(P,\mu_P)$ is a projective $A^e$-resolution  of $A$, $\Delta_P$ is a comultiplication on it, and $M$ is an $A$-algebra, then we can define a map 
$$\smile_{\Delta_P}:\Hom_{A^e}(P,M)\otimes\Hom_{A^e}(P,M)\rightarrow\Hom_{A^e}(P,M)$$ by the equality $f\smile_{\Delta_P}g=\mu_M(f\otimes g)\Delta_P$. 
This multiplication induces a well defined product on $\HH^*(A,M)$, which coincides with the cup product $\smile$. Note also that the bar resolution 
of $A$ 
admits a comultiplication $\Delta$ defined by the equality
$$
\Delta(a_0\otimes\cdots\otimes a_{n+1})=\sum\limits_{i=0}^n(a_0\otimes\cdots a_i\otimes 1_A)\otimes_A(1_A\otimes a_{i+1}\otimes\cdots\otimes a_{n+1}).
$$\\
Note also that the multiplication $\smile_{\Delta}$ on $C(A,M)$ coincides with $\smile$ and that 
$$(\Delta\otimes\Id_{\Bar(A)})\Delta=(\Id_{\Bar(A)}\otimes\Delta)\Delta.$$

\section{Group Cohomology}\label{Ho}

In this section we review some facts on group cohomology. Most of the facts stated here can be generalized for the setting of Hopf algebras.

Let $G$ be a group. For $M,N\in\Mod\kk G$ we define the $\kk G$-module structure  on $M\otimes N$ and $\Homk(M,N)$ by the formulas
 $$
 \alpha(x\otimes y)=(\alpha x)\otimes (\alpha y);\,\,(\alpha f)(x)=\alpha\big(f(\alpha^{-1}x)\big) \label{notation}
 $$
 for $x\in M$, $y\in N$, $f\in\Homk(M,N)$ and $\alpha\in G$. For $\alpha\in G$ we denote by $p_{\alpha}\in \kk G^*$ the map such that 
$p_{\alpha}(\alpha)=1$ and $p_{\alpha}(\beta)=0$ for $\beta\in G\setminus\{\alpha\}$.

We define next a special $\kk G$-projective resolution of the trivial module $\kk$. Which we will denote by $\Bar\kk$. Since $\kk$ can be considered
as an algebra, there are two distint objects, which are denoted by $\Bar\kk$, nevertheless there will be no 
confusion, since, in this paper, $\Bar\kk$ will always have the meaning below. 

\begin{Def}[\bf Bar Resolution of the Trivial Module, and Group Cohomology]
 
 The bar resolution $(\Bar(\kk), \sum_{\alpha\in G}p_{\alpha})$ of the trivial $\kk G$-module $\kk$ is the projective resolution  defined in the 
following way. Its $n$-th 
component is $\Bar_n(\kk)=\kk G^{\otimes(n+1)}$. The differential $d_{\Bar(\kk)}$ is defined by the equality
\begin{multline*}
d_{\Bar(\kk)}(\alpha_0\otimes \alpha_1\otimes\cdots\otimes \alpha_n)=\sum\limits_{i=0}^{n-1}(-1)^i\alpha_0\otimes\cdots\otimes \alpha_{i-1}\otimes 
\alpha_i\alpha_{i+1}\otimes \alpha_{i+2}\otimes\cdots\otimes \alpha_{n}\\
+(-1)^n\alpha_0\otimes \alpha_1\otimes\cdots\otimes \alpha_{n-1}
\end{multline*}
for $n>0$ and $\alpha_i\in G$ ($0\le i\le n$). The cohomology of the group $G$ with coefficients in the $\kk G$-module $M$ is the holomogy of the complex 
$C(G,M)=\Hom_{\kk G}(\Bar(\kk), M)$.
\end{Def}

Note that $\Bar(\kk)=\Bar(\kk G)\otimes_{\kk G} \kk$. We write $C^n(G,M)$ instead of $C_{-n}(G,M)$ and $\delta^n_M$ instead of $d_{C(G,M),-1-n}$. Then $C^n(G,M)\cong 
\Hom_{\kk}(\kk G^{\otimes n}, M)$ and $\delta^n:C^n(G,M)\rightarrow C^{n+1}(G,M)$ is defined by the equality
\begin{multline*}
\delta^n_M(\hat{u})(\alpha_1\otimes\cdots\otimes \alpha_{n+1})=\alpha_1\hat{u}(\alpha_2\otimes\cdots\otimes \alpha_{n+1})\\
+\sum\limits_{i=1}^n(-1)^i\hat{u}(\alpha_1\otimes\cdots\otimes \alpha_{i-1}\otimes \alpha_i\alpha_{i+1}\otimes \alpha_{i+2}\otimes\cdots\otimes 
\alpha_{n+1})+(-1)^{n+1}\hat{u}(\alpha_1\otimes\cdots\otimes \alpha_n)
\end{multline*}
for $\hat{u}\in C^n(G,M)$, for $\alpha_i\in G$ ($1\le i\le n+1$). By definition the group cohomogy of $G$ with coefficients in $M$ is $\Ho^*(G,M)=\oplus_{n\ge 
0}\Ho^n(G,M)$, where $\Ho^n(G,M)=\Ker\delta^n_M/\Im\delta^{n-1}_M$.  It follows from our definition that $\Ho^n(G, M)= \Ext^n_{\kk G}(\kk, M).$

Let $A$ be some $\kk$-algebra and $G$ be a group. We say that $G$ acts on $A$ if we fix some group homomorphism $\eta:G\rightarrow \Aut_{\kk}(A)$.
For an action of a group $G$ on an algebra $A$, we write ${}^{\alpha}\!a$ instead of $\eta(\alpha)(a)$ for $\alpha\in G$, $a\in A$. Note that any vector space can be equipped by trivial multiplication and so we define what is an action of a group on a vector space.
If a group $G$ acts on some vector space $V$,  then $V^G$ denotes the subspace $\{v\in V\mid {}^{\alpha}\!v=v\,\,\forall \alpha\in G\}$ of $V$.

If $G$ acts on an algebra $M$, then there is a product $\smile:C(G,M)\otimes C(G,M)\rightarrow C(G,M)$, 
where for $\hat{u}\in C^n(G,M)$ and 
$\hat{v}\in C^m(G,M)$ the element $\hat{u}\smile\hat{v}\in C^{n+m}(G)$ is defined by the equality
$(\hat{u}\smile\hat{v})(\alpha_1\otimes\cdots\otimes \alpha_{n+m})=\hat{u}(\alpha_1\otimes\cdots\otimes 
\alpha_{n})\big(\alpha_1\dots\alpha_n\hat{v}(\alpha_{n+1}\otimes\cdots\otimes \alpha_{n+m})\big)$. This product determines a well-defined product 
$\smile$ on $\Ho^*(G,M)$

Let $(Q,\mu_Q)$ be a $\kk G$-projective resolution of $\kk$ and $M$ be a $\kk G$-module, a morphism $\Delta_Q:Q\rightarrow Q\otimes Q$ of $\kk 
G$-complexes will be called  a comultiplication on $Q$ if it 
satisfies the 
equality $\mu_{\kk}(\mu_Q\otimes\mu_Q)\Delta_Q=\mu_Q$.

Assume that $G$ acts on an algebra $M$. Then we can define a map $$\smile_{\Delta_Q}:\Hom_{\kk 
G}(Q,M)\otimes\Hom_{\kk 
G}(Q,M)\rightarrow\Hom_{\kk G}(Q,M)$$ by the equality $f\smile_{\Delta_Q}g=\mu_M(f\otimes g)\Delta_Q$. This map induces a well defined product on 
$\Ho^*(G,M)$, which coincides with the product $\smile$. Note that if we have a nonassociative multiplication on $M$ (but such that the multiplication 
by 
$\alpha$ is an automorphism for any $\alpha\in G$), then we still can define the map $\smile_{\Delta_Q}$ on $\Hom_{\kk G}(Q,M)$ and $\Ho^*(G,M)$. In 
this case the obtained algebra will not be associative, but the algebra structure on $\Ho^*(G,M)$ still does not depend on the choice of $Q$ and 
$\Delta_Q$.

Note that the bar resolution of $\kk$ admits a comultiplication $\Delta$ defined by the equality
$$
\Delta(\alpha_0\otimes \alpha_1\otimes\cdots\otimes \alpha_n)=\oplus_{i=0}^n(\alpha_0\otimes\cdots \otimes\alpha_i)\otimes(\alpha_0\cdots 
\alpha_i\otimes \alpha_{i+1}\otimes\cdots\otimes \alpha_n).
$$
It is clear that the map $\Delta$ induces a map $\smile:C(G,M)\otimes C(G,M)\rightarrow C(G,M)$.

Note that analogously we can define an $\Ho^*(G,\kk)$-module structure on $\Ho^*(G,M)$ for any $\kk G$-module $M$. By the Golod-Venkov-Evens theorem 
(see \cite{Golod}, \cite{Venkov} and \cite{Evens}) if $|G|<\infty$, then the algebra $\Ho^*(G,\kk)$ is finitely generated and for any finitely 
generated $\kk G$-module $M$ the $\Ho^*(G,\kk)$-module $\Ho^*(G,M)$ is finitely generated.

\section{Smash product}\label{smash}

In this section we recall some results on the smash product. In this paper there are two situations where it appears, which are, as actions of groups associated with the group algebra and as $G$-gradings associated with $\kk G^*$.

Let $H$ be a bialgebra with the counit $\ee:H\rightarrow \kk$ and the comultiplication $\Delta:H\rightarrow H\otimes H$. We use the standard notation for the
comultiplication $\Delta(h)=\sum_{(h)}h_{(1)}\otimes h_{(2)}$.

\begin{Def}The unital algebra $A$ is an $H$-module algebra if there is a map $\cdot:H\otimes A\rightarrow A$ such that\\
(1) $\cdot$ determines an $H$-module structure on $A$,\\
(2) $h\cdot ab=\sum_{(h)}(h_{(1)}\cdot a)(h_{(2)}\cdot b)$ ($a,b\in A$, $h\in H$),\\
(3) $h\cdot 1_A=\ee(h)1_A$ ($h\in H$).
\end{Def}

\begin{Def}
Let $A$ be an $H$-module algebra. We define the smash product $A\#H$ in the following way. The algebra $A\#H$ is isomorphic to $A\otimes H$ as a vector 
space. The multiplication of the elements $a\otimes \alpha, b\otimes \beta\in A\#H$ is defined by the formula
$$
(a\otimes \alpha)(b\otimes \beta)=\sum\limits_{(\alpha)}a(\alpha_{(1)}\cdot b)\otimes(\alpha_{(2)}\beta).
$$
From now on we write $a\alpha$ instead of $a\otimes \alpha$.
\end{Def}

\begin{Ex}
Let $H=\kk G$ for some group $G$. In this case an $H$-module algebra structure on $A$ corresponds to a $G$-action. We write simply $AG$ instead of 
$A\#\kk G$. Then $\Delta(\alpha)=\alpha\otimes \alpha$ and the multiplication on $AG$ is defined by the formula
$$
(a\alpha)(b\beta)=(a{}^{\alpha}\!b)(\alpha\beta)\,\,(a,b\in A, \alpha,\beta\in G).
$$
\end{Ex}
We say that $A$ is $G$-graded if there is some direct sum decomposition $A=\oplus_{\alpha\in 
G}A_{\alpha}$ such that $A_{\alpha}A_{\beta}\subset A_{\alpha\beta}$ for all $\alpha,\beta\in G$. In this case $a_{\alpha}$ ($a\in A$, $\alpha\in G$) 
denotes the image of $a$ under the canonical projection from $A$ to $A_{\alpha}$.

\begin{Ex}\label{covering}
Let $H=\kk G^*$ for some finite group $G$. In this case an $H$-module algebra structure on $B$ corresponds to a $G$-grading. We write simply $BG^*$ 
instead of $B\#\kk G^*$. The elements $p_{\alpha}$ ($\alpha\in G$) constitute a basis of $\kk G^*$. Moreover, $\Delta(p_{\alpha})=\sum_{\beta\in 
G}p_{\alpha\beta^{-1}}\otimes p_{\beta}$ and the multiplication on $BG^*$ is defined by the formula
$$
(ap_{\alpha})(bp_{\beta})=(ab_{\alpha\beta^{-1}}) p_{\beta}\,\,(a,b\in B, \alpha,\beta\in G).
$$
\end{Ex}

\begin{Def}[{\bf Graded Bimodule, Homomorphism of Graded Bimodule}]
Let $B$ be a $G$-graded algebra. We say that a $B$-bimodule $M$ is $G$-graded if $M$ is equipped with a direct sum decomposition $M=\oplus_{\alpha\in 
G}M_{\alpha}$ such that $B_{\alpha} M_{\beta} B_{\gamma} \subset M_{\alpha\beta\gamma}$ for all triples $(\alpha,\beta,\gamma) \in G\times G\times G$..

A homomorphism of $G$-graded bimodules is a homomorphism of bimodules that preserves $G$-grading.
\end{Def}

 A  resolution  $(P,\mu_P)$  of $B$ is 
called a $G$-graded projective resolution if each bimodule $P_i$ is equipped with a $G$-grading in such a way that the differential $d_P$ and the map 
$\mu_P$ are morphisms of $G$-graded bimodules. For example, $\Bar(B)$ becomes a $G$-graded resolution of $B$ if we define 
$(B^{\otimes(n+2)})_\alpha$ as the 
space generated by such $a_0\otimes\cdots\otimes a_{n+1}$ that $a_i\in B_{\alpha_i}$ and $\alpha_0\dots\alpha_{n+1}=\alpha$. \\

\begin{Def}
Let $B$ be a $G$-graded algebra and  $M$ and $N$ be $G$-graded $B$-bimodules, we define $\Hom_{B^e}(M,N)_0$ as the vector subspace of $\Hom_{B^e}(M,N)$ consisting of
grading preserving homomorphisms, i.e.\\
$\{f \in \Hom_{B^e}(M,N) \mid f(x)_{\alpha}=f(x)$ for any $\alpha\in G$ and any $x\in M_{\alpha}$\}.\\
We also define a complementary subspace $\Hom_{B^e}(M,N)_1$ as the set\\
$\{f \in \Hom_{B^e}(M,N) \mid f(x)_{\alpha}=0$ for any $\alpha\in G$ and any $x\in M_{\alpha} \}$.
\end{Def}

The following equality holds:
$\Hom_{B^e}(M,N) =\Hom_{B^e}(M,N)_0 \oplus \Hom_{B^e}(M,N)_1$.

Let us fix  $G$-graded projective bimodule resolution $(P,\mu_P)$ of $B$.
Let $\Hom_{B^e}(P,B)_{k}$ denote $\oplus_{i\ge 0}\Hom_{B^e}(P_i,B)_k$ for $k=0,1$. Then it is easy to see that there is a direct sum decomposition 
$$\Hom_{B^e}(P,B)=\Hom_{B^e}(P,B)_0\oplus \Hom_{B^e}(P,B)_1$$
of complexes. It defines a decomposition $\HH^*(B)=\HH^*(B)_0\oplus \HH^*(B)_1$. It is easy to see that this decomposition does not depend on the 
choice of $G$-graded resolution.

Moreover, it is clear that we can construct a comultiplication $\Delta_P$ on $P$ that at the same time is a morphism of $G$-graded bimodules.
Then the algebra $\Hom_{B^e}(P,B)_0$ turns into a subalgebra of $\Hom_{B^e}(P,B)$ and $\HH^*(B)_0$ turns into a subalgebra of $\HH^*(B)$. Moreover, 
$\Hom_{B^e}(P,B)_1$ is a $\Hom_{B^e}(P,B)_0$-bimodule and $\HH^*(B)_1$ is a $\HH^*(B)_0$-module.

Let $A=BG^*$. In this case there is a natural $G$-action on $A$  defined by the equality ${}^{\alpha}\!(ap_{\beta})=ap_{\beta\alpha^{-1}}$ ($a\in 
A$, 
$\alpha,\beta\in G$). The algebra $AG$ has a canonical $G$-grading defined by the equality $(AG)_{\alpha}=A\alpha$ for $\alpha\in G$.
It is well known that $AG$ is Morita equivalent to $B$ (see \cite{CohMon}). This equivalence can be defined in the following way. Firstly, by results 
of \cite{CohMon} we have $A\cong p_1AGp_1$. It is easy to see that it turns into isomorphism of graded algebras if we define 
$(p_1AGp_1)_{\alpha}=p_1A\alpha p_1$ for $\alpha\in G$ (see the proof of \cite[Lemma 3.4]{CohMon}). Then the $AG$-bimodule $M$ corresponds to 
$p_1AGp_1$-bimodule $p_1Mp_1$ and a morphism of $AG$-bimodules $f:M\rightarrow N$ corresponds to the morphism $p_1fp_1:p_1Mp_1\rightarrow p_1Np_1$ of 
$p_1AGp_1$-bimodules that sends $x\in p_1Mp_1$ to $f(x)\in p_1Np_1$. It is clear that the constructed functor gives a functor between the categories of 
$G$-graded bimodules and so the isomorphism of graded algebras $\HH^*(B)\cong\HH^*(AG)$ induces isomorphisms $\HH^*(B)_0\cong\HH^*(AG)_0$ and 
$\HH^*(B)_1\cong\HH^*(AG)_1$.\\

The Morita equivalence  between the algebras above were generalized in the setting of  $\kk$-categories with group actions, in \cite{CM}.

\section{$G$-action on Hochschild complex}\label{GActH}

In this section we work with some fixed algebra $A$ equipped with an action of the group $G$. In particular, $A$ and $A^e$ are $\kk G$-modules. Note 
also that $\Hom_A(M,N)$ is a $\kk G$-submodule of $\Homk(M,N)$ if $M,N\in\Mod AG$. 

There is an inclusion of algebras $\phi: A^eG\rightarrow (AG)^e$ given by  $\phi\big((a\otimes 
b)\alpha\big)=a\alpha\otimes\alpha^{-1}b$ for $a,b\in A$, $\alpha\in G$.
If $M$ is an $AG$-bimodule, then we also write ${}^{\alpha}\!x$ instead of $\alpha x\alpha^{-1}$ for 
$\alpha\in G$, $x\in M$.
In this case, the inclusion above defines by restriction of scalar an $A^eG$-module structure on $M$, which is given
by the formula $(a\otimes b)\alpha x=a{}^{\alpha}\!xb$ for $a,b\in A$, $\alpha\in G$ and $x\in 
M$. 
Note also that $A$ is a left $A^eG$-module with multiplication defined by the formula $(a\otimes b)\alpha c=a{}^{\alpha}\!cb$ for $a,b,c\in A$ and 
$\alpha\in G$. If $M$ is an $A^eG$-module and $N$ has a structure of $\kk G$-module, then we define the $A^eG$-module structure on $M\otimes N$ by the 
equality $a\alpha(x\otimes y)=(a\alpha x)\otimes (\alpha y)$ for $a\in A^e$, $\alpha\in G$, $x\in M$, $y\in N$. At the same time, if $M$ is a left 
$AG$-module and $N$ is a right $AG$-module, then for any vector space $V$ we define the $(AG)^e$-module structure on $M\otimes V\otimes N$ by the 
equality $(a\otimes b)(x\otimes v\otimes y)=(ax)\otimes v\otimes (yb)$ for $a,b\in AG$, $x\in M$, $y\in N$ and $v\in V$.
\begin{Def}[\bf $A^eG$-module Algebra]
 An $A^eG$-module $M$ will be 
called an $A^eG$-module algebra if $M$ is a $\kk G$-module algebra and $A$-algebra at the same time.
\end{Def}
 
 Let  $M$ be in $\Mod A^eG$. Then $C^n(A,M)$ is equipped with a $\kk G$-module structure described above. 

Let us recall now some results of \cite{Volk5}.
 
Firstly, the $\kk G$-module structure on $C^n(A,M)$ induces a $G$-action on $\HH^*(A,M)$.
In addition, if there is an algebra structure on $M$ such that it becomes an $A^eG$-module algebra, then we obtain $G$-action on the graded algebras 
$(C(A,M),\smile)$ and $(\HH^*(A,M),\smile)$.
 The $G$-action on $C(A,M)$ defines new complex $C(A,M)^G$ with differential $\delta_M^G:C(A,M)^G\rightarrow C(A,M)^G$, whose homology we 
denote by $\HH^n(A,M)^{G\uparrow}=(\Ker\delta_M^{G,n})/(\Im\delta_M^{G,n-1})$. If $M$ is an $A^eG$-module algebra, then $\HH^*(A,M)^{G\uparrow}$ is a graded 
$\kk$-algebra and there is a homomorphism of algebras $\Theta_{M}^G:\HH^*(A,M)^{G\uparrow}\rightarrow\HH^*(A,M)^{G}$ induced by the inclusion of 
$C(A,M)^G\hookrightarrow C(A,M)$. Finally, if $G$ is a finite group such that $\charr\kk\nmid|G|$, then $\Theta_{M}^G$ is an isomorphism.

 Let now $(P,\mu_P)$ be an $A^e$-projective $A^eG$-resolution of $A$. Then $\Hom_{A^e}(P,M)$ is a $\kk G$-complex, in this case we get a 
$G$-action on $\HH^n(A,M)$  which coincides with the action defined above by \cite[Lemma 2]{Volk5}.

Consider next the complex $\Hom_{A^eG}(P,M)$. Let 
$$\HH^n(A,M)^{G\uparrow}_{P}=\Ker\big(\Hom_{A^eG}(d_{P,n},M)\big)/\Im\big(\Hom_{A^eG}(d_{P,n-1},M)\big)$$
denote its homology. The inclusion $\Hom_{A^eG}(P,M)\hookrightarrow \Hom_{A^e}(P,M)$ induces a map 
$\Theta_{P,M}^G:\HH^*(A,M)^{G\uparrow}_P\rightarrow\HH^*(A,M)^{G}$.
 In particular, it is clear that $\HH^*(A,M)^{G\uparrow}_{\Bar(A)}=\HH^*(A,M)^{G\uparrow}$. Note also that for $m\ge 0$ we can apply the functor $\Ho^m(G,-)$ to the 
$\kk G$-complex $\Hom_{A^e}(P,M)$. We denote by $\HH^n(A,M)^{G,m}_P$ the $n$-th homology of the complex $\Ho^m(G, \Hom_{A^e}(P,M))$. It is clear that 
$\HH^n(A,M)^{G,0}_P=\HH^n(A,M)^{G\uparrow}_P$, since $\Ho^0(G, \Hom_{A^e}(P,M))= \Hom_{A^eG}(P,M)$.

Assume now that $M$ is an $A^eG$-module algebra. Assume additionally that there is a comultiplication $\Delta_P:P\rightarrow P\otimes_A P$, which is a 
morphism of $A^eG$-modules. Then we have a map $\smile_{\Delta_P}:\Hom_{A^eG}(P,M)\otimes \Hom_{A^eG}(P,M)\rightarrow \Hom_{A^eG}(P,M)$. If the maps 
$(\Delta_P\otimes\Id_P)\Delta_P$ and $(\Id_P\otimes\Delta_P)\Delta_P$ are homotopic as morphisms of $A^eG$-complexes, then the map $\smile_{\Delta_P}$ 
induces a well-defined multiplication $\smile_{\Delta_P}$ on $\HH^*(A,M)^{G\uparrow}_{P}$. It is clear that if $P$ is an $A^eG$-projective resolution 
of $A$, then the map $\Delta_P$ exists and induces a well-defined product on $\HH^*(A,M)^{G\uparrow}_{P}$. Moreover, in this case any comultiplication 
on $P$, which is a morphism of  $A^eG$-modules, defines the same product on $\HH^*(A,M)^{G\uparrow}_{P}$.  We can see  that for any 
$A^e$-projective $A^eG$-resolution 
$(P,\mu_P)$ of $A$ and any comultiplication $\Delta_P$ the map $\Theta_{P,M}^G$ is a homomorphism of graded algebras (one of which can be 
nonassociative in general).

Usually the role of $M$ is played by $A$ or $AG$. It is well known that $\HH^*(A)$ is a graded commutative algebra. The next result shows that this is 
not true for $\HH^*(A,AG)$ in general. On the other hand, we will show that the algebra $\HH^*(A,AG)^G$ is graded commutative.
Note that there is a direct sum decomposition $\HH^*(A,AG)=\oplus_{\alpha\in G}\HH^*(A,A\alpha)$.

\begin{lemma}
Let $\bar{f}$ be an element of $\HH^i(A,A\alpha)$ and $\bar{g}$ be an element of $\HH^j(A,A\beta)$, where $i,j\ge 0$, $\alpha,\beta\in G$. Then the 
equality
$$
\bar{f}\smile \bar{g}=(-1)^{ij}\bar{g}\smile{}^{\beta^{-1}}\!\bar{f}
$$
is satisfied in $\HH^*(A,AG)$.
\end{lemma}
\begin{Proof} The proof, in fact, is analogous to  the proof of \cite[Theorem 3]{Gerstenhaber}, we give it here for the sake of completeness. 

Let $f\in C^i(A,A\alpha)$ and 
$g\in C^j(A,A\beta)$ be the elements representing $\bar{f}$ and $\bar{g}$ correspondingly.
Note that $g(X)\beta^{-1}$ lies in $A$ for any $X\in A^{\otimes n}$ and so we can define the elements $f\circ_n g\in C^{i+j-1}(A,A\alpha\beta)$ ($1\le 
n\le i$) by the equality
\begin{multline*}
(f\circ_n g)(a_1\otimes\cdots\otimes a_{i+j-1})\\
=(-1)^{n(j-1)}f(a_1\otimes\cdots \otimes a_{n-1}\otimes g(a_n\otimes\cdots\otimes a_{n+j-1})\beta^{-1}\otimes {}^{\beta}\!a_{n+j}\otimes\cdots\otimes 
{}^{\beta}\!a_{i+j-1})\beta
\end{multline*}
for $a_k\in A$ ($1\le k\le i+j-1$).
Let us define $f\circ g=\sum\limits_{n=1}^if\circ_ng.$
Note that
$$
\begin{aligned}
 \delta_{i+j-1}^k(f\circ_n g)&=\begin{cases}
(-1)^{j-1}\delta_i^k(f)\circ_{n+1} g,&\mbox{if $0\le k\le n-1$},\\
(-1)^{n-1}f\circ_{n} \delta_j^{k-n+1}(g),&\mbox{if $n\le k\le n+j-1$},\\
(-1)^{j-1}\delta_i^{k-j+1}(f)\circ_n g,&\mbox{if $n+j\le k\le i+j$},
 \end{cases}\\
\delta_i^{n-1}(f)\circ_n g&=\begin{cases}
(-1)^{j+n}f\circ_{n-1} \delta_j^0(g),&\mbox{if $2\le n\le i+1$},\\
(-1)^{j-1}g\smile{}^{\beta^{-1}}\!f,&\mbox{if $n=1$},
 \end{cases}\\
 \delta_i^n(f)\circ_n g&=\begin{cases}
(-1)^{j+n+1}f\circ_n \delta_j^{j+1}(g),&\mbox{if $1\le n\le i$}.\\
(-1)^{(i+1)j}f\smile g,&\mbox{if $n=i+1$},
 \end{cases}\\
 \end{aligned}$$
 Let now apply $\delta_{AG}^{i+j-1}$ to the element $f\circ g$:
\begin{multline*}
\delta_{AG}^{i+j-1}(f\circ g)=\sum\limits_{k=0}^{i+j}\sum\limits_{n=1}^i\delta_{i+j-1}^k(f\circ_n 
g)=\sum\limits_{n=1}^i\big(\sum\limits_{k=0}^{n-1}(-1)^{j-1}\delta_i^k(f)\circ_{n+1} g\\
+\sum\limits_{k=n}^{n+j-1}(-1)^{n-1}f\circ_{n} \delta_j^{k-n+1}(g)+\sum\limits_{k=n+j}^{i+j}(-1)^{j-1}\delta_i^{k-j+1}(f)\circ_n g\big)\\
=(-1)^{j-1}\sum\limits_{n=1}^{i+1}(\delta_{AG}^{i}(f)-\delta_i^{n-1}(f)-\delta_i^n(f))\circ_n g\\
+\sum\limits_{n=1}^i(-1)^{n-1}f\circ_n(\delta_{AG}^{j}(g)-\delta_j^0(g)-\delta_j^{j+1}(g))=-g\smile{}^{\beta^{-1}}\!f+(-1)^{ij}f\smile g.
\end{multline*}
Thus, the lemma is proved.
\end{Proof}

\begin{coro}
$\HH^*(A,AG)^G$ is a graded commutative algebra.
\end{coro}

\section{From $A^eG$ to $(AG)^e$}

In this section we explain how the $(AG)^e$-projective resolution of $AG$ can be obtained from $A^eG$-projective resolution of $A$. 
Also it is noted above that there is a homomorphism of algebras $\phi: A^eG\rightarrow (AG)^e$ defined by the formula $\phi\big((a\otimes 
b)\alpha\big)=a\alpha\otimes\alpha^{-1}b$ for $a,b\in A$, $\alpha\in G$. This inclusion makes $(AG)^e$ an $A^eG$-projective module and induces a 
forgetful functor from $\Mod (AG)^e$ to $\Mod 
A^eG$, for which there is a left adjoint functor ${\bf T}=(AG)^e\otimes_{A^eG}-$. So we have natural isomorphisms 
$$\chi_{N,M}:\Hom_{A^eG}(N,M)\rightarrow\Hom_{(AG)^e}({\bf T}N,M)$$ for all $M\in \Mod (AG)^e$ and $N\in \Mod A^eG$.

\begin{lemma}\label{AeGtoAGe}
Let $(P,\mu_P)$ be an $A^eG$-projective resolution of $A$. Then ${\bf T}P$ is an $(AG)^e$-projective resolution of $AG$. In particular, the map 
$\chi_{P,M}$ induces a bijective map $\bar\chi_{P,M}:\HH^*(A,M)^{G\uparrow}_{P}\rightarrow \HH^*(AG,M)$ for any $AG$-bimodule $M$. Moreover, if $M$ is 
an $AG$-algebra, then $\bar\chi_{P,M}$ is an isomorphism of graded algebras.
\end{lemma}
\begin{Proof}  It is clear that the functor ${\bf T}$ is exact and sends projective $A^eG$-modules to projective $(AG)^e$-modules. So it remains to 
prove the last assertion. Let $\Delta_P:P\rightarrow P\otimes_A P$ be a comultiplication, which is a morphism of $A^eG$-modules.
Let $\psi:{\bf T}(P\otimes_A P)\rightarrow {\bf T}P\otimes_{AG} {\bf T}P$ be defined by the equality
$$\psi\big((a\otimes b)\otimes_{A^eG}(x\otimes_A y)\big)=\big((a\otimes 1_A)\otimes_{A^eG} x\big)\otimes_{AG} \big((1_A\otimes b)\otimes_{A^eG} 
y\big)$$
for $a,b\in AG$ and $x,y\in P$.
 Then we can define a comultiplication $\Delta_{{\bf T}P}$ by the equality $\Delta_{{\bf T}P}=\psi{\bf T}\Delta_P$. Note that for $x\in P$, $a\in 
(AG)^e$ and $f\in \Hom_{A^eG}(P,M)$ we have $\chi_{P,M}(f)(a\otimes x)=af(x)$.
 
 Let us prove that $$\chi_{P,M}\big(\mu_M(f\otimes g)\big)=\mu_M\big(\chi_{P,M}(f)\otimes\chi_{P,M}(g)\big)\psi:{\bf T}(P\otimes_A P)\rightarrow M$$ 
for $f,g\in\Hom_{A^eG}(P,M)$. For $a,b\in AG$ and $x,y\in P$ we have
 $$
 \chi_{P,M}\big(\mu_M(f\otimes g)\big)\big((a\otimes b)\otimes(x\otimes y)\big)=(a\otimes b)f(x)g(y)=af(x)g(y)b.
 $$
 and
\begin{multline*}
 \mu_M\big(\chi_{P,M}(f)\otimes\chi_{P,M}(g)\big)\psi\big((a\otimes b)\otimes(x\otimes y)\big)\\
 =\mu_M\big(\chi_{P,M}(f)\otimes\chi_{P,M}(g)\big)\Big(\big((a\otimes 1_A)\otimes x\big)\otimes \big((1_A\otimes b)\otimes y\big)\Big)=af(x)g(y)b.
\end{multline*}
Then we have
\begin{multline*}
\chi_{P,M}\big(\mu_M(f\otimes g)\Delta_P\big)=\chi_{P,M}\big(\mu_M(f\otimes g)\big){\bf T}\Delta_P\\
=\mu_M\big(\chi_{P,M}(f)\otimes\chi_{P,M}(g)\big)\psi{\bf T}\Delta_{P}=\mu_M\big(\chi_{P,M}(f)\otimes\chi_{P,M}(g)\big)\Delta_{{\bf T}P}.
\end{multline*}
\end{Proof}

\begin{Def} We call the $A^eG$-module $P$ almost $A^eG$-projective if the module $P\otimes \kk G$ is $A^eG$-projective.
\end{Def}

It is clear that a module $P$ in $\Mod A^eG$ is almost $A^eG$-projective iff $P\otimes Q$ is $A^eG$-projective for any $\kk G$-projective module $Q$. 
Note that for 
any $\kk G$-module $M$ the module $A^e\otimes M$ is almost $A^eG$-projective. In particular,  the bar resolution of $A$ is almost $A^eG$-projective, 
because $A^{\otimes(n+2)}\cong A^e\otimes A^{\otimes n}$ as $A^eG$-module.

The next lemma can be found in \cite[Proposition 1.7]{ARS}.
\begin{lemma}\label{proj_alg}
Let $G$ be a finite group and $A$ be a $\kk G$-algebra. Then the following statements are equivalent:\\
{\rm 1. }$A$ is $\kk G$-projective;\\
{\rm 2. }there is such $a\in A$ that $\sum_{\alpha\in G}{}^{\alpha}a=1_A$;\\
{\rm 3. }$A$ is $AG$-projective.
\end{lemma}

\begin{coro}\label{proj}
Let $G$ be a finite group and $A$ be a $\kk G$-algebra. If $A$ is a $\kk G$-projective module, then the bar resolution of $A$ is an $A^eG$-projective 
resolution.
\end{coro}
\begin{Proof} Since $A$ is $\kk G$-projective, $A^e$ is $\kk G$-projective too. So $A^e$ is a projective $A^eG$-module by Lemma \ref{proj_alg}. So 
$\Bar_n(A)\cong A^e\otimes A^{\otimes n}$ is $A^eG$-projective.
\end{Proof}

If $|G|<\infty$ and $A$ is a projective $\kk G$-module, then it follows from Lemma \ref{AeGtoAGe} and Corollary \ref{proj} that there is an isomorphism 
of graded spaces $\chi_{\Bar(A),M}:\HH^*(A,M)^{G\uparrow}\rightarrow \HH^*(AG,M)$ for any $AG$-bimodule $M$. Moreover, $\chi_{\Bar(A),M}$ is  an 
algebra isomorphism if $M$ is an $AG$-algebra.

\section{Constructing spectral sequences from the zero page}\label{multspec}

Let $(P, \mu_P)$  be an almost $A^eG$-projective resolution of $A$ and $(Q,\mu_Q)$ be a $\kk G$-projective resolution of the trivial $\kk G$-module $\kk$. 
Then $P\otimes Q$ is an $A^eG$-projective resolution of $A$ by the arguments above. We can compute the Hochschild cohomology of $AG$ with coefficients 
in an $(AG)^e$-module $M$ using the complex $\Hom_{A^eG}(P\otimes Q,M)$ accordingly to Lemma \ref{AeGtoAGe}.

\begin{Def} A double $A$-complex is a $\mathbb{Z}\times\mathbb{Z}$-graded $A$-module $P$ with two differentials of degrees $(-1,0)$ and $(0,-1)$, i.e. 
an $A$-module $P$ with some fixed $A$-module direct sum decomposition $P=\oplus_{(n,m)\in \mathbb{Z}\times \mathbb{Z}}P_{n,m}$ and $A$-module 
homomorphisms $d_P^h,d_P^v:P\rightarrow P$ such that $d_P^h(P_{n,m})\subset P_{n-1,m}$, $d_P^v(P_{n,m})\subset P_{n,m-1}$ and 
$(d_P^h)^2=(d_P^v)^2=d_P^hd_P^v+d_P^vd_P^h=0$. The totalization of $P$ is the complex $\Tot(P)$, which equals $P$ as $A$-module, has grading 
$\Tot(P)_n=\oplus_{i+j=n}P_{i,j}$ and differential $d_P^h+d_P^v$.
\end{Def}

As a tensor product of complexes $P$ and $Q$ the complex $P\otimes Q$ is a totalization of the double complex $B$, whose $(i,j)$-member is 
$B_{i,j}=P_i\otimes Q_j$ and differentials are defined by the equalities $d_B^h(x\otimes y)=d_P(x)\otimes y$ and $d_B^v(x\otimes y)=(-1)^ix\otimes 
d_Q(y)$ for $x\in P_i$, $y\in Q$. Aplying the functor $\Hom_{A^eG}(-,M)$ to the double complex $B$ we obtain the double complex $\Hom_{A^eG}(B,M)$ 
whose totalization is isomorphic to $\Hom_{A^eG}(P\otimes Q,M)$.  As any double complex $\Hom_{A^eG}(B,M)$ has two filtrations: column filtration and 
row filtration. Correspondingly we obtain two spectral sequences. Let ${}^I\!E$ and ${}^{II}\!E$ denote these spectral sequences. We write $\delta_M^h$ and $\delta_M^v$ instead of $\Hom_{A^eG}(d_B^h,M)$ and $\Hom_{A^eG}(d_B^v,M)$.

Let firstly consider the sequence ${}^I\!E$ that is obtained from the filtration
$$
\cdots\subset F_1^I\subset F_0^I= \Hom_{A^eG}(B,M),\,\,F_n^I=\Hom_{A^eG}(\oplus_{i\ge n, j\ge 0}B_{i,j},M).
$$
The zero page of ${}^I\!E$ is defined by the equality ${}^I\!E_0^{i,j}=\Hom_{A^eG}(P_j\otimes Q_i,M)$. There are natural isomorphisms 
$$\phi_{i,j}:\Hom_{A^eG}(P_j\otimes Q_i,M)\cong \Hom_{\kk G}(Q_i,\Hom_{A^e}(P_j,M))$$ defined by the equality
$\big((\phi_{i,j} f)(x)\big)(y)=f(y\otimes x)$ for $f\in \Hom_{A^eG}(P_j\otimes Q_i,M)$, $x\in Q_i$, $y\in P_j$.  Because of that one sees that the first 
page is ${}^I\!E_1^{i,j}\cong \Hom_{\kk G}(Q_i,\HH^j(A,M))$ and the second page is
${}^I\!E_2^{i,j}\cong \Ho^i(G,\HH^j(A,M))$. Since spectral sequence of double complex converges to the homology of its totalization, we obtain that
\begin{equation}\label{spec_1}
{}^I\!E_2^{i,j}\cong \Ho^i(G,\HH^j(A,M))\implies \HH^{i+j}(AG,M).
\end{equation}
This fact is well known, its analogues can be found in \cite{Sanada}, \cite{Stefan} or \cite{CibRed}. Note that the authors of \cite{CibRed} use 
similar constructions  to the ones we make here.

Let now consider the sequence ${}^{II}\!E$ that is obtained from the filtration
$$
\cdots\subset F_1^{II}\subset F_0^{II}= \Hom_{A^eG}(B,M),\,\,F_n^{II}=\Hom_{A^eG}(\oplus_{i\ge 0, j\ge n}B_{i,j},M).
$$
The zero page of ${}^{II}\!E$ is defined by the equality
$${}^{II}\!E_0^{i,j}=\Hom_{A^eG}(P_i\otimes Q_j,M)\cong \Hom_{\kk G}(Q_j,\Hom_{A^e}(P_i,M)).$$
It is clear that ${}^{II}\!E_1^{i,j}\cong \Ho^j(G,\Hom_{A^e}(P_i,M))$ and ${}^{II}\!E_2^{i,j}\cong \HH^i(A,M)_P^{G,j}$. So
\begin{equation}\label{spec_2}
{}^{II}\!E_2^{i,j}\cong \HH^i(A,M)_P^{G,j}\implies \HH^{i+j}(AG,M).
\end{equation}
If $\Hom_{A^e}(P,M)$ is a $\kk G$-projective module, then ${}^{II}E_2^{i,j}=0$ for $j>0$ and this means that any element of 
$\HH^n(A,M)^{G\uparrow}_{P\otimes Q}\cong\HH^n(AG,M)$ can be represented by an element of $\Hom_{A^eG}((P\otimes Q)_n,M)$, which equals 0 on 
$P_i\otimes Q_j$ for $j>0$. In particular, we obtain  that $\Theta_{P,M}^G:\HH^*(A,M)^{G\uparrow}_P\rightarrow\HH^*(A,M)^{G}$ is an isomorphism if 
$\Hom_{A^e}(P,M)$ is a $\kk G$-projective module.
Note that $\Hom_{A^e}(P,M)$ can be $\kk G$-projective even if $P$ is not $A^eG$-projective. For example, it is so if $P$ is the bar resolution of $A$ and $M$ is 
a $\kk G$-projective module.

It is easy to see that the composition of edge maps 
$${}^{II}\!E_2^{*,0}\twoheadrightarrow{}^{II}\!E_{\infty}^{*,0}\hookrightarrow\HH^*(AG,M)\twoheadrightarrow{}^I\!E_{\infty}^{0,*}\hookrightarrow 
{}^I\!E_{2}^{0,*}$$ is simply $\Theta_{P,M}^G$. If $\Hom_{A^e}(P,M)$ is $\kk G$-projective, then 
${}^{II}\!E_2^{*,0}\twoheadrightarrow{}^{II}\!E_{\infty}^{*,0}\hookrightarrow\HH^*(AG,M)$ is an isomorphism and so ${}^I\!E_{\infty}^{0,*}\cong 
\Im\Theta_{P,M}^G$.

Let ${}^I\!d_r:{}^I\!E_r\rightarrow {}^I\!E_r$ and ${}^{II}\!d_r:{}^{II}\!E_r\rightarrow {}^{II}\!E_r$ denote the $r$-th differentials of spectral 
sequences.

 \section{Multiplication on the zero page}\label{multsec}
 
Let now $M$ be an $AG$-algebra. Let $P$ and $Q$ be as in the previous section. Assume additionally that we have some comultiplication $\Delta_P$ on 
$P$, which is a morphism of $A^eG$-modules, and some comultiplication $\Delta_Q$ on $Q$. It always exists, because $Q$ is a $\kk G$-projective 
resolution. Let $$\psi: (P\otimes_A P)\otimes (Q\otimes Q)\rightarrow (P\otimes Q)\otimes_A (P\otimes Q)$$ be defined by the equality 
$\psi(x_1\otimes_A x_2\otimes y_1\otimes y_2)=(-1)^{ij}x_1\otimes y_1\otimes_A x_2\otimes y_2$ for $x_1\in P$, $x_2\in P_i$, $y_1\in Q_j$, $y_2\in Q$. 
It is easy to see that $\psi$ is an isomorphism of $A^eG$-complexes. Let us define $\Delta_{P\otimes Q}=\psi(\Delta_P\otimes\Delta_Q)$. Then 
$\Delta_{P\otimes Q}$ is a comultiplication on $P\otimes Q$, which is a morphism of $A^eG$-complexes.
Thus, we can define the multiplication $\smile_{\Delta_{P\otimes Q}}$ on the complex $\Hom_{A^eG}(P\otimes Q,M)$. Since $P\otimes Q$ is 
$A^eG$-projective, we can describe the multiplication $\smile$ on the Hochschild cohomology of $AG$ describing the multiplication induced by 
$\smile_{\Delta_{P\otimes Q}}$ on $\HH^*(A,M)_{P\otimes Q}^{G\uparrow}$  accordingly to Lemma \ref{AeGtoAGe}.

On the other hand $\Delta_{P\otimes Q}$ induces an $A^eG$-homomorphism from $B$ to $B\otimes_A B$, which respects the 
$\mathbb{Z}\times\mathbb{Z}$-grading and  both differentials of $B$. So $\smile_{\Delta_{P\otimes Q}}$ induces a map
$$\smile_B:\Hom_{A^eG}(B,M)\otimes \Hom_{A^eG}(B,M)\rightarrow \Hom_{A^eG}(B,M),$$
which respects the $\mathbb{Z}\times\mathbb{Z}$-grading and the both differentials of $\Hom_{A^eG}(B,M)$. Such a map induces maps 
$\smile^I_r:{}^I\!E_r\otimes {}^I\!E_r\rightarrow {}^I\!E_r$ and $\smile^{II}_r:{}^{II}\!E_r\otimes {}^{II}\!E_r\rightarrow {}^{II}\!E_r$
for all $0\le r\le \infty$. Thus, we obtain two spectral sequences of algebras $({}^I\!E_r,\smile^I_r)$ and $({}^{II}\!E_r,\smile^{II}_r)$, and both of 
them converge to $ \HH^*(AG,M)$.

Note that modulo isomorphism $\Hom_{A^eG}(B,M)\cong \Hom_{\kk G}(Q, \Hom_{A^e}(P,M))$ the multiplication $\smile_B$ is simply the multiplication 
$\smile_{\Delta_Q}$ on the complex $\Hom_{\kk G}(Q, L)$, where $L$ is the $\kk G$-module algebra $\Hom_{A^e}(P,M)$ with the multiplication $\smile_{\Delta_P}$. If $\Delta_P$ and 
$\Delta_Q$ are coassociative, then $\smile_{\Delta_P}$, $\smile_{\Delta_Q}$ and $\smile_B$ are associative, but in general all of this multiplications may be 
nonassociative. On the other hand, the multiplication induced by $\smile_B$ on ${}^I\!E_2$ is the multiplication $\smile$ on the complex $\Ho^*(G, \bar 
L)$, where $\bar L$ is the $\kk G$-module algebra $\HH^*(A,M)$ with the multiplication $\smile$. So $\smile_B$ always induces an associative multiplication on 
${}^I\!E_2$.

\begin{rema}
Let $P'$ be some $A^eG$-projective resolution of $A$ and $Q'$ be some $\kk G$-projective resolution of $\kk$. Let ${}^I\!E'$ and ${}^{II}\!E'$ be the spectral 
sequences constructed using $P'$ and $Q'$. There are a homomorphism of $A^eG$-complexes $\Phi:P'\rightarrow P$ lifting $\Id_A$ and a homomorphism of 
$\kk G$-complexes $\Psi:Q'\rightarrow Q$ lifting $\Id_{\kk}$. Then the map $\Phi\otimes\Psi:P'\otimes Q'\rightarrow P\otimes Q$ induces some maps 
${}^I\Upsilon_r:{}^I\!E_r\rightarrow{}^I\!E_r'$ and ${}^{II}\Upsilon_r:{}^{II}\!E_r\rightarrow{}^{II}\!E_r'$ for all $r\ge 0$. It is easy to see that 
the map ${}^I\Upsilon_2$ is an isomorphism of algebras (for any choice of $\Delta_{P'}$ and $\Delta_{Q'}$). So the map ${}^I\Upsilon_r$ is an 
isomorphism of algebras for any $r\ge 2$ (see \cite[Theorem 3.4]{guide}). This means that beginning from the second page the sequence ${}^I\!E$ does 
not depend on the choice of $P$ and $Q$.
\end{rema}

\begin{Ex}
Let us take $P=\Bar(A)$, $\Delta_P=\Delta$, $Q=\Bar(\kk)$ and $\Delta_Q=\Delta$  (see the definitions in Sections \ref{HH} and \ref{Ho}). Then
$$\Hom_{A^eG}(B_{i,j},M)\cong C^j\big(G,C^i(A,M)\big)\cong C^i(A,M)\otimes C^j(G,\kk).$$
For $f\in C^i(A,M)$, $\hat{u}\in C^j(G,\kk)$ we have
$$
\delta^h_M(f\otimes \hat{u})=\delta_M^i(f)\otimes \hat{u},\,\,\delta^v_M(f\otimes \hat{u})=(-1)^i\big(f\otimes\delta^j(\hat{u})+\sum\limits_{\alpha\in 
G}({}^{\alpha}\!f-f)\otimes(p_{\alpha}\smile\hat{u})\big).
$$
Finally, the map $\smile_B$ is defined by the equality
$$
(f\otimes \hat{u})\smile_B(g\otimes \hat{v})=(-1)^{j_1i_2}\sum_{\alpha\in G}(f\smile{}^{\alpha}\!g)\otimes(\hat{u}_{\alpha}\smile \hat{v})
$$
for $f\in C^{i_1}(A,M)$, $g\in C^{i_2}(A,M)$, $\hat{u}\in C^{j_1}(G,\kk)$, and $\hat{v}\in C^{j_2}(G,\kk)$, where
$$\hat{u}_{\alpha}(\alpha_1\otimes\dots\otimes\alpha_{j_1})=\begin{cases}\hat{u}(\alpha_1\otimes\dots\otimes\alpha_{j_1}),&\mbox{if 
$\alpha_1\dots\alpha_{j_1}=\alpha$},\\0,&\mbox{otherwise.}\end{cases}$$
\end{Ex}

For a finite subgroup $H$ of $G$ such that $\charr\kk\nmid|H|$ we define $e_H=\frac{\sum_{\alpha\in H}\alpha}{|H|}$. If $H\unlhd G$, then $e_H$ is a
central idempotent and the subalgebra $(\kk G)e_H$ is isomorphic to $\kk (G/H)$. Moreover, if $(Q,\mu_Q)$ is a projective $\kk G$-resolution of $\kk$, 
then $(e_HQ,\mu_Q|_{e_HQ_0})$ is a projective $\kk G$-resolution of $\kk$ too. Assume that $Q=e_HQ$ initially. Then $Q$ is a $\kk(G/H)$-resolution and 
we have
\begin{equation}\label{eH}
\Hom_{A^eG}(B,M)\cong \Hom_{\kk G}(Q, \Hom_{A^e}(P,M))\cong \Hom_{\kk(G/H)}(Q,e_H\Hom_{A^e}(P,M)).\end{equation}

\begin{Ex}\label{cyclic}
In this example we take $P=\Bar(A)$ and $\Delta_P=\Delta$ again.
Assume that $\charr \kk=p$, where $p$ is some prime. Let $q$ be a power of $p$ and $C_q$ be a cyclic group of order $q$. Let $G$ be an extension of 
$C_q$ by a finite group $H$, whose order is not divisible by $p$, i.e. $H\unlhd G$ and $G/H\cong C_q$. Let $\rho$ be a generator of $C_q$. Let us 
denote by $\NN$ the element $\sum\limits_{i=0}^{q-1}\rho^i$. Let us define a $\kk G$-projective resolution $(Q,\sum_{\alpha\in G}p_{\alpha}|_{(\kk 
G)e_H})$ of $\kk$. We define $Q_n=(\kk G)e_H$ for all $n\ge 0$. Given some $x\in(\kk G)e_H$, we write $x[n]$ for the corresponding element of $Q_n$. 
Then the differential $d_Q$ is defined by the equality
$$
d_Q(e_H[n])=\begin{cases}
\NN e_H[n-1],&\mbox{if $2\mid j$ and $j>0$,}\\
(e_H-\rho e_H)[n-1],&\mbox{if $2\nmid j$.}
\end{cases}
$$
Accordingly to \eqref{eH} we have $\Hom_{A^eG}(B_{i,j},M)\cong C^i(A,M)^H$. Given some $f\in C^i(A,M)^H$, we write $f[j]$ for the corresponding element 
of $\Hom_{A^eG}(B_{i,j},M)$.
Then for $f\in C^i(A,M)^H$ we have
$$
\delta^h_M(f[j])=\delta_M^i(f)[j],\,\,
\delta^v_M(f[j])=\begin{cases}
(-1)^i(f-{}^{\rho}\!f)[j+1],&\mbox{if $2\mid j$,}\\
(-1)^i\sum\limits_{t=0}^{q-1}{}^{\rho^t}\!f[j+1],&\mbox{if $2\nmid j$.}
\end{cases}
$$
Let now define the map $\Delta_Q:Q\rightarrow Q\otimes Q$ by the equalities
\begin{multline*}
\Delta_Q(e_H[2n])=\sum\limits_{m=0}^ne_H[2m]\otimes e_H[2(n-m)]\\
+\sum\limits_{m=0}^{n-1}\sum\limits_{i=0}^{q-1}\sum\limits_{j=0}^{i-1}\rho^ie_H[2m+1]\otimes \rho^je_H[2(n-m)-1],\\
\Delta_Q(e_H[2n+1])=\sum\limits_{m=0}^ne_H[2m]\otimes e_H[2(n-m)+1]+\sum\limits_{m=0}^n\rho e_H[2m+1]\otimes e_H[2(n-m)]
\end{multline*}
for $n\ge 0$. Then the map $\smile_B$ is defined by the equality
$$
f[j_1]\smile_Bg[j_2]=\begin{cases}
(f\smile g)[j_1+j_2],&\mbox{if $2\mid j_1$,}\\
(-1)^{i_2}({}^{\rho}\!f\smile g)[j_1+j_2],&\mbox{if $2\nmid j_1$ and $2\mid j_2$,}\\
(-1)^{i_2}\left(\sum\limits_{n=0}^{q-1}\sum\limits_{m=0}^{n-1}{}^{\rho^n}\!f\smile{}^{\rho^m}\!g\right)[j_1+j_2],&\mbox{if $2\nmid j_1$ and $2\nmid j_2$.}\\
\end{cases}
$$
for $f\in C^{i_1}(A,M)^H$, $g\in C^{i_2}(A,M)^H$.
\end{Ex}

\section{When ${}^I\!E$ degenerates}

\begin{Def}
Spectral sequence $E$ is called $(R,S)$-degenerated if $E_{R}^{i,j}=0$ for $i\ge S$. In particular, $(R,S)$-degenerated spectral sequence 
collapses at the $max(R,S)$ page.
\end{Def}

Let $G$ be a finite group and $M$ be an $A^eG$-module algebra. In this section we consider the following question. Do integers $R$ and $S$ such that 
${}^I\!E$ is $(R,S)$-degenerated exist?

Let us formulate an enough condition for $(R,S)$-degenerating of ${}^I\!E$. We call the family of elements $\hat{u}_t\in\Ho^{i_t}(G,\kk)$ ($t\in \XX$) 
an $S_0$-generating system if $\Ho^i(G,L)=\sum\limits_{t\in\XX}\Ho^{i-i_t}(G,L)\smile\hat{u}_t$ for any $i\ge S_0$ and any simple $\kk G$-module $L$. 
Since there are only finitely many simple $\kk G$-modules, then by the Golod-Venkov-Evens theorem there is a finite $S_0$-generating system for some 
integer $S_0$.

\begin{theorem}\label{degeneration}
Let $\hat{u}_t\in\Ho^{i_t}(G,\kk)$ ($1\le t\le s$) be an $S_0$-generating system for $\Ho^*(G,\kk)$. Assume that $R\ge 2$ is some integer. Then\\
{\rm 1.} if the sequence ${}^I\!E$ is $(R,S)$-degenerated, then the class of $\hat{u}_t^{k_t}\smile 1_{\Ho^*(G,\HH^*(A,M))}$ in ${}^I\!E_{R}^{i_tk_t,0}$ 
equals $0$ for any $1\le t\le s$, where $k_t=\upper{\frac{S}{i_t}}$;\\
{\rm 2.} if the integers $k_t$ ($1\le t\le s$) are such that the class of $\hat{u}_t^{k_t}1_{\Ho^*(G,\HH^*(A,M))}$ in ${}^I\!E_{R}^{i_tk_t,0}$ equals 
$0$ for any $1\le t\le s$, then ${}^I\!E$ is $(R,S)$-degenerated for some $S\le S_0+\sum\limits_{t=1}^s(k_t-1)i_t$.
\end{theorem}
\begin{Proof} The first assertion is obvious. We write simply $\hat{u}_t$ instead of $\hat{u}_t\smile 1_{\Ho^*(G,\HH^*(A,M))}$. Let us prove that ${}^I\!E_R^{i,j}=\sum\limits_{t=1}^s{}^I\!E_R^{i-i_t,j}\smile_I^R\hat{u}_t$ for $i\ge 
S_0$. We may assume that $(Q,\mu_Q)$ is a minimal $\kk G$-projective resolution of $\kk$.
Let $u_t:Q_{i_t}\rightarrow \kk$ ($1\le t\le s$) be an element representing $\hat{u}_t$. There are morphisms of complexes $T_t:Q\rightarrow Q[-i_t]$ 
such that $\mu_QT_{t,i_t}=u_t$. Let us define $U=\oplus_{t=1}^sQ[-i_t]$. Let $\Phi=\sum\limits_{t=1}^sT_t:Q\rightarrow U$. Let us take some $i\ge S_0$. 
Let $J$ be the Jacobson radical of $\kk G$ and $\pi_i:Q_i\rightarrow Q_i/JQ_i$ be the canonical projection. Since $Q$ is a minimal resolution of $\kk$, 
the map $\pi_i$ defines a cocycle $\hat{\pi}_i\in\Ho(G,Q_i/JQ_i)$.
Since $i\ge S_0$, it is clear that $\Ho^i(G,L)=\sum\limits_{t=1}^s\Ho^{i-i_t}(G,L)\smile\hat{u}_t$ for any semisimple $\kk G$-module $L$. Hence, there 
are some $\hat{v}_t\in \Ho^{i-i_t}(G,Q_i/JQ_i)$ such that $\hat{\pi}_i=\sum\limits_{t=1}^s\hat{v}_t\smile\hat{u}_t$. Let $v_t:Q_{i-i_t}\rightarrow \kk$ 
($1\le t\le s$) be an element representing $\hat{v}_t$ and let $\Gamma_i=\sum\limits_{t=1}^sv_t:U_i\rightarrow Q_i/JQ_i$. It is easy to see that 
$\Gamma_i\Phi_i$ represents $\hat{\pi}_i$ in $\Ho(G,Q_i/JQ_i)$. Since $Q_i/JQ_i$ is semisimple and $Q$ is minimal, we have $\pi_i=\Gamma_i\Phi_i$. 
Since $U$ is projective, there is a map $\bar\Psi_i:U_i\rightarrow Q_i$ such that $\pi_i\bar\Psi_i=\Gamma_i$, i.e. 
$\pi_i(\Id_{Q_i}-\bar\Psi_i\Phi_i)=0$. Then the image of $f_i=\Id_{Q_i}-\bar\Psi_i\Phi_i$ lies in $JQ_i$ and so $\bar\Psi_i\Phi_i$ is invertible. Let 
us define $\Psi_i=(\bar\Psi_i\Phi_i)^{-1}\bar\Psi_i$. We have $\Psi_i\Phi_i=\Id_{Q_i}$. It is easy to show that $\Psi_id_{U,i}=d_{Q,i}\Psi_{i+1}$ for 
$i\ge S_0$. Let now $\bar{w}$ be an element of ${}^I\!E_R^{i,j}$, where $i\ge S_0$. Then there is some element $w\in\Hom_{\kk G}(Q_i, \Hom_{A^e}(P_j,M))$ that 
represents $\bar{w}$. We have $w=(w\Psi_i)\Phi_i=\sum\limits_{t=1}^s(w\Psi_i)|_{Q_{i-i_t}}T_{t,i}$. Then it is not hard to show that 
$(w\Psi_i)|_{Q_{i-i_t}}$ represents some element $\hat{w}_t\in {}^I\!E_R^{i-i_t,j}$ for any $1\le t\le s$ and 
$\bar{w}=\sum\limits_{t=1}^s\hat{w}_t\smile_I^R\hat{u}_t$ in ${}^I\!E_R^{i,j}$.

 Assume now that the condition of the second assertion is satisfied.
Let $w$ be an element of ${}^I\!E_R^{i,j}$, where $i\ge S_0+\sum\limits_{t=1}^s(k_t-1)i_t$. It is easy to show that 
$w=\sum\limits_{t=1}^sv_t\smile_I^R\hat{u}_t^{k_t}$ for some $v_t\in{}^I\!E_R^{i-i_tk_t,j}$. Sinse $\hat{u}_t^{k_t}=0$ in ${}^I\!E_R^{i_tk_t,0}$,
$w=0$ in ${}^I\!E_R^{i,j}$.
\end{Proof}

\begin{rema}
Note that it follows from the proof of Theorem \ref{degeneration} that for any $\kk G$-module $L$ the $\Ho^*(G,\kk)$-module $\Ho^*(G,L)$ is generated 
by elements whose degrees are less or equal than $S_0$, where $S_0$ is a minimal integer such that $\HH^*(G,\kk)$ is an $S_0$-generating system.
\end{rema}

\begin{coro}\label{period_deg}
Let $R\ge 2$ be some integer. If $Q_s$ is periodic with period $n$, then ${}^I\!E$ is $(R,S)$-degenerated for some $S$ iff 
$\hat{u}^m1_{\Ho^*(G,\HH^*(A,M))}=0$ in ${}^I\!E_R^{mn,0}$ for some $m$, where $\hat{u}$ is an element of $\Ho^n(G,\kk)$ represented by the map 
$\mu_Q:Q_n\rightarrow \kk$. If $m$ is the minimal integer satisfying the required equality, then the minimal $S$ such that ${}^I\!E$ is 
$(R,S)$-degenerated is more than $(m-1)n$ and less or equal than $mn$.
\end{coro}
\begin{Proof} It is easy to see that $\hat{u}$ is an $n$-generating system. So everything follows from Theorem \ref{degeneration}.
\end{Proof}

\begin{coro} If there are such $T$ and $R\ge 2$ that ${}^I\!E_R^{i,0}=0$ for $i\ge T$, then ${}^I\!E$ is $(R,S)$-degenerated for some integer $S$. If there 
is such $T$ that ${}^I\!E_{\infty}^{i,0}=0$ for $i\ge T$, then ${}^I\!E$ is $(R,S)$-degenerated for some integers $R$ and $S$
\end{coro}

\section{Some facts about ${}^I\!E$ in the case of cyclic group}

From this moment to the end of the paper we assume that $\charr\kk=p$ and $G$ is a group from Example \ref{cyclic}. This means that there is a finite 
normal subgroup $H$ of $G$ such that $p\nmid|H|$ and $G/H\cong C_q$, where $q$ is a power of $p$. As in Example \ref{cyclic} we denote by $\rho$ the 
generator of $C_q$. Let $\Delta_{\rho}=1-\rho\in\kk C_q$. Note that $\sum\limits_{i=0}^{q-1}\rho^i=\Delta_{\rho}^{q-1}$. In this section we formulate 
and prove some properties of ${}^I\!E$, which follow from the results of previous sections.

\begin{prop}\label{simple1} Let $A$ be a $\kk G$-module algebra, $M$ be an $A^eG$-module algebra. Then the following statements are true:\\
{\rm 1.} If ${}^I\!E$ is $(R,S)$-degenerated for some integers $R$ and $S$, then $M^H$ is $\kk C_q$-projective.\\
{\rm 2.} Let $M^{AH}$ be a $\kk C_q$-nonprojective module. Then ${}^I\!E_{\infty}^{1,0}\not=0$. In particular, we have $\dim_{\kk}\HH^1(AG,M)>0$ in the 
case where $M$ is an $(AG)^e$-module.\\
{\rm 3.} Let $M^{AH}$ be a $\kk C_q$-projective module. Then ${}^I\!E$ is $(2,1)$-degenerated. In particular, $\HH^*(AG, M)\cong\HH^*(A,M)^G$ as graded 
algebra in the case where $M$ is an $(AG)^e$-module.\\
\end{prop}
\begin{Proof} 1. We may assume that the sequence ${}^I\!E$ is constructed using $P$ and $Q$ from Example \ref{cyclic}. Also we use the notation 
introduced there. There is some integer $k$ such that the image of $\Ho^{2k}(C_q,M^{AH})$ in ${}^I\!E_R^{2k,0}$ is zero. In particular, the class of 
$1_M\in M^{AG}\subset M^{AH}\cong\Hom_{\kk C_q}(Q_{2k},M^{AH})$ in ${}^I\!E_R^{2k,0}$ is zero. Let us prove that for any $0\le i\le 2k$ there is such 
$u_i\in M^{H}$ that the class of $1_M-\Delta_{\rho}^{q-1}u_i$ in ${}^I\!E_{i+1}^{2k,0}$ is zero. We can take $u_{2k}=0$. Suppose that we have 
constructed $u_{i+1}$. Then there is some $\bar{h}_i\in {}^I\!E_{i+1}^{2k-i-1,i}$ satisfying the equality 
${}^I\!d_{i+1}(\bar{h}_i)=1_M-\Delta_{\rho}^{q-1}u_{i+1}$ in ${}^I\!E_{i+1}^{2k,0}$. Let $\bar{h}_i$ be represented by $h_i\in C^i(A,M)^H$. Then there 
are such $h_j\in C^j(A,M)^H$ ($0\le j\le i-1$) that $\delta_M^h(h_j[2k-1-j])=\delta_M^v(h_{j+1}[2k-2-j])$ for $0\le j\le i-1$. Then 
${}^I\!d_{i+1}(\bar{h}_i)$ is represented by $\delta_M^v(h_0[2k-1])=\Delta_{\rho}^{q-1}h_0[2k]$. So we can take $u_i=u_{i+1}+h_0$. Now we know that 
$1_M-\Delta_{\rho}^{q-1}u_0$ equals $0$ in ${}^I\!E_{1}^{2k,0}=\Hom_{\kk}(Q_{2k},M^{AH})\cong M^{AH}$, i.e. $1_M=\Delta_{\rho}^{q-1}u_0$ for some 
$u_0\in M^{H}$. By Lemma \ref{proj_alg} the $\kk C_q$-module $M^H$ is projective.

2. The statement is obvious, since ${}^I\!E_{\infty}^{1,0}={}^I\!E_2^{1,0}=\Ho^1(C_q,M^{AH})$ and $\Ho^1(C_q,L)\not=0$ for any $\kk C_q$-nonprojective 
module $L$. Note that this statement is true for any $A^eG$-module $M$.

3. It is enough to prove that $\HH^*(A,M)^H$ is a $\kk C_q$-projective module. By Lemma \ref{proj_alg} there is an element $a\in M^{AH}\subset 
\HH^*(A,M)^H$ such that $\sum\limits_{\alpha\in C_q}{}^{\alpha}a=1_{M^{AH}}=1_{\HH^*(A,M)^H}$. Using Lemma \ref{proj_alg} again we obtain that 
$\HH^*(A,M)^H$ is $\kk C_q$-projective.
\end{Proof}

Now we prove a theorem that gives some connections between the dimensions of $\HH^n(AG,M)$ and the dimensions of $\HH^n(A,M)^G$.

\begin{theorem}\label{dims}
Let $A$ be a $\kk G$-module algebra, $M$ be an $(AG)^e$-module. Then $$\dim_{\kk}\HH^n(AG,M)\ge\dim_{\kk}\HH^n(A,M)^G$$ for any $n\ge 0$. Moreover, if 
$M$ is additionally an $A^eG$-module algebra, then the following conditions are equivalent:\\
{\rm 1.} ${}^I\!E$ is $(3,2)$-degenerated;\\
{\rm 2.} $\dim_{\kk}\HH^n(AG,M)=\dim_{\kk}\HH^n(A,M)^G$ for any $n\ge 0$;\\
{\rm 3.} $\dim_{\kk}\HH^1(AG,M)=\dim_{\kk}\HH^1(A,M)^G$.\\
\end{theorem}
\begin{Proof} Let $(P,\mu_P)$ be an $A^eG$-projective resolution of $A$. Let us introduce some notation for this proof. Let $C^n$ and $\delta^n$ denote 
$\Hom_{A^e}(P_n,M)$ and $d_{\Hom(A^e,P),-1-n}$ correspondingly. Note that $C^n$ is a $\kk G$-projective module. Let $Z^n=\Ker\delta^n\subset C^n$ and 
$B^n=\Im\delta^n\subset C^{n+1}$. During this proof we write simply $\HH^n$ instead of $\HH^n(A,M)$. Note also that for any $\kk C^q$-module $L$ and 
any $n\ge 1$ the equality $\dim_{\kk}\Ho^n(C_q,L)=\dim_{\kk}\Ho^1(C_q,L)$ holds, since both parts of the equality equal $0$ for projective indecomposible 
$L$ and equal $1$ for nonprojective indecomposible $L$. So 
$$\dim_{\kk}\Ho^n(G,L)=\dim_{\kk}\Ho^n(C_q,L^H)=\dim_{\kk}\Ho^1(C_q,L^H)=\dim_{\kk}\Ho^1(G,L)$$ for any $\kk G$-module $L$ and any $n\ge 1$.

We have the following two exact sequences of $\kk G$-modules:
\begin{equation}\label{exact}
Z^{n-1}\rightarrowtail C^{n-1}\twoheadrightarrow B^{n-1}\mbox{ and }B^{n-1}\rightarrowtail Z^n\twoheadrightarrow \HH^n.
\end{equation}
From the second  of them we obtain the long exact sequence
\begin{equation}\label{long_exact}
(B^{n-1})^G\rightarrowtail (Z^n)^G\rightarrow (\HH^n)^G\rightarrow\Ho^1(G, B^{n-1})\rightarrow\cdots
\end{equation}
From the first sequence from \eqref{exact} we obtaint an exact sequence
\begin{equation}\label{4_exact}
(Z^{n-1})^G\rightarrowtail (C^{n-1})^G\rightarrow (B^{n-1})^G\twoheadrightarrow\Ho^1(G, Z^{n-1})
\end{equation}
 and an isomorphism $\Ho^1(G,B^{n-1})\cong \Ho^2(G,Z^{n-1})$. Therefore, we obtain from \eqref{long_exact} that
\begin{multline*}
 \dim_{\kk}(\HH^n)^G+\dim_{\kk}(B^{n-1})^G\le  \dim_{\kk}(Z^n)^G+\dim_{\kk}\Ho^1(G, B^{n-1})\\
 =\dim_{\kk}(Z^n)^G+\dim_{\kk}\Ho^2(G, Z^{n-1})=\dim_{\kk}(Z^n)^G+\dim_{\kk}\Ho^1(G, Z^{n-1}).
\end{multline*}
 On the other hand it follows from \eqref{4_exact} that
 $$
\dim_{\kk}(Z^{n-1})^G+\dim_{\kk}(B^{n-1})^G= \dim_{\kk}(C^{n-1})^G+\dim_{\kk}\Ho^1(G, Z^{n-1}).
 $$
 So we have
 $$
 \dim_{\kk}(\HH^n)^G\le  \dim_{\kk}(Z^n)^G+\dim_{\kk}(Z^{n-1})^G-\dim_{\kk}(C^{n-1})^G.
 $$
 We know from Lemma \ref{AeGtoAGe} that $\HH^n(AG,M)\cong \HH^n(A,M)^{G\uparrow}_{P}$. By definition of $\HH^n(A,M)^{G\uparrow}_{P}$ we have an exact 
sequence
 $$
(Z^{n-1})^G\rightarrowtail (C^{n-1})^G\rightarrow (Z^n)^G\twoheadrightarrow \HH^n(A,M)^{G\uparrow}_{P}.
 $$
 So
 $$
 \dim_{\kk}\HH^n(AG,M)=\dim_{\kk}(Z^n)^G+\dim_{\kk}(Z^{n-1})^G-\dim_{\kk}(C^{n-1})^G\ge  \dim_{\kk}\HH^n(A,M)^G.
 $$
 
 Let now prove the second part of the theorem. Note that the implication $"2\Rightarrow 3"$ is obvious.
 
 $"1\Rightarrow 2."$ For any $n\ge 2$ we have two exact sequences:
$$
{}^I\!E_3^{n,0}\rightarrowtail (\HH^n)^G\rightarrow \Ho^2(G,\HH^{n-1})\rightarrow\cdots\rightarrow \Ho^{2n-2}(G,\HH^1)\twoheadrightarrow 
\Ho^{2n}(G,\HH^0)
$$
 and
$$
 {}^I\!E_3^{n-1,1}\rightarrowtail \Ho^1(G,\HH^{n-1})\rightarrow \Ho^3(G,\HH^{n-2})\rightarrow\cdots\rightarrow \Ho^{2n-3}(G,\HH^1)\twoheadrightarrow 
\Ho^{2n-1}(G,\HH^0).
$$
So we have
\begin{multline*}
\dim_{\kk}\HH^n(AG,M)=\dim_{\kk}{}^I\!E_3^{n,0}+\dim_{\kk}{}^I\!E_3^{n-1,1}\\
=\dim_{\kk}(\HH^n)^G+\sum\limits_{i=1}^n(-1)^i\dim_{\kk}\Ho^{2i}(G,\HH^{n-i})+\sum\limits_{i=1}^n(-1)^{i-1}\dim_{\kk}\Ho^{2i-1}(G,\HH^{n-i})\\
=\dim_{\kk}(\HH^n)^G+\sum\limits_{i=1}^n(-1)^i\big(\dim_{\kk}\Ho^{2i}(G,\HH^{n-i})-\dim_{\kk}\Ho^{2i-1}(G,\HH^{n-i})\big)\\
=\dim_{\kk}\HH^n(A,M)^G.
\end{multline*}
 
 $"3\Rightarrow 1."$ We have 
\begin{multline*}
\dim_{\kk}\HH^1(AG,M)= \dim_{\kk}\Ho^1(G,\HH^0)+\dim_{\kk}{}^I\!E_3^{0,1}\\
=\dim_{\kk}\Ho^1(G,\HH^0)+\dim_{\kk}(\HH^1)^G+\dim_{\kk}{}^I\!E_3^{2,0}-\dim_{\kk}\Ho^2(G,\HH^0)\\
=\dim_{\kk}\HH^1(A,M)^G+\dim_{\kk}{}^I\!E_3^{2,0}.
\end{multline*}
So we have ${}^I\!E_3^{2,0}=0$. Hence, ${}^I\!E$ is $(3,2)$-degenerated by Corollary \ref{period_deg}.
\end{Proof}

\begin{rema} Note that Theorem \ref{dims} becomes true for any $A^eG$-module $M$ if we change everywhere $\HH^n(AG,M)$ by $\HH^n(A,M)^{G\uparrow}_{P}$, where $(P,\mu_P)$ is some $A^eG$-projective resolution of $A$.
\end{rema}

\section{(3,2)-degeneration of ${}^I\!E$}

Now we are going to explore the case where the sequence ${}^I\!E$ is $(3,2)$-degenerated. Let us recall that for an algebra $A$ and an $A$-bimodule $M$ 
a $\kk$-linear map $f$ from $A$ to $M$ is called derivation if $f(ab)=af(b)+f(a)b$ for all $a,b\in A$. Let $\Der(A,M)$ denote the set of all 
derivations from $A$ to $M$. For $x\in M$ we denote by $\phi_x$ the inner derivation defined by $x$, i.e. $\phi_x:A\rightarrow M$ is the map defined by 
the equality $\phi_x(a)=ax-xa$ for $a\in A$.

\begin{lemma}\label{3_2_deg}
Let $A$ be a $\kk G$-module algebra, $M$ be an $A^eG$-module algebra. Then the sequence ${}^I\!E$ is $(3,2)$-degenerated iff there are such $\mm\in M^H$ 
and $\hh\in\Der(A, M)^H$ that
 $\Delta_{\rho}^{q-1}\mm=1_M$ and $\Delta_{\rho}\hh+\phi_{\mm}=0$.
\end{lemma}
\begin{Proof} We may assume that the sequence ${}^I\!E$ is constructed using $P$ and $Q$ from Example \ref{cyclic}. By Corollary \ref{period_deg} the 
sequence ${}^I\!E$ is $(3,2)$-degenerated iff the class of $1_M\in M^{AG}\cong\Hom_{\kk G}(Q_2,M^A)$ in ${}^I\!E^{2,0}_3$ equals $0$. Let ${}^I\!E$ be 
$(3,2)$-degenerated. This means that there is such element $\bh\in\HH^1(A,M)^G$ that ${}^I\!d_2(\bh)$ equals to the class of $1_M$ in 
$\Ho^2(G,M^{AH})$. Let $\hh\in C^1(A,M)^H$ be a map representing $\bh$. Since $\delta_M^1(\hh)=0$, $\hh$ is an element of $\Der(A, M)^H$. Since 
$\bh\in\HH^1(A,M)^G$, there is some element $\mm_0\in C^0(A,M)^H=M^{H}$ satisfying the equality $\phi_{\mm_0}=\Delta_{\rho}\hh$. Then the class of 
${}^I\!d_2(\bh)$ equals to the class of $-\Delta_{\rho}^{q-1}\mm_0$ in $\Ho^2(G,M^{AH})$. So there is some element $\mm_1\in M^{AH}$ satisfying the 
equality $\Delta_{\rho}^{q-1}(\mm_1-\mm_0)=1_M$. Let us take $\mm=\mm_1-\mm_0$. Since $\phi_{\mm_1}=0$, the required $\mm$ and $\hh$ are constructed. Conversely, if 
the elements $\mm$ and $\hh$ satisfy the conditions of the theorem, then $\hh$ represents some element $\bh\in\HH^1(A,M)^G$ such that ${}^I\!d_2(\bh)$ equals 
to the class of $1_M$ in $\Ho^2(G,M^{AH})$ and so ${}^I\!E$ is $(3,2)$-degenerated.
\end{Proof}

Let $M$ be an $A^eG$-module algebra. If $\mm\in M^H$ and $\hh\in\Der(A, M)^H$ satisfy the conditions of Lemma \ref{3_2_deg}, then we call $M$ an 
$(\mm,\hh)$-degenerated algebra. Note that any $(\mm,\hh)$-degenerated algebra is a $\kk G$-projective module. In particular, the map $\chi_{\Bar(A),M}$ is an 
isomorphism in this case and ${}^I\!E_{\infty}^{0,*}\cong\Im\Theta_M^G$ by arguments of Section \ref{multspec}.

\begin{rema}\label{change}
If  $M$ is an $(\mm,\hh)$-degenerated algebra and $\mm'\in M^H$ is such that $\Delta_{\rho}^{q-1}(\mm')=1_M$, then $M$ is an $(\mm',\hh')$-degenerated algebra for 
some $\hh'$. Really, since $M^H$ is $\kk C_q$-projective module and $\Delta_{\rho}^{q-1}(\mm-\mm')=0$, there is such $x\in M^H$ that $\Delta_{\rho}x=\mm-\mm'$ 
and we can take $\hh'=\hh+\phi_x$.
\end{rema}

Let $M$ be an $(\mm,\hh)$-degenerated algebra. Then we have the following commutative diagram of algebras with exact rows:
$$
\xymatrix{
\Ker\Theta_M^G\,\ar@{>->}[r]&\HH^*(A,M)^{G\uparrow}\ar@{=}[d]\ar@{->>}[r]&\Im\Theta_M^G\ar@{->}[d]^{\cong}\\
{}^I\!E_{\infty}^{1,*-1}\,\ar@{>->}[r]&\HH^*(A,M)^{G\uparrow}\ar@{->>}[r]&{}^I\!E_{\infty}^{0,*}.
}
$$
This means that $\Ker\Theta_M^G={}^I\!E_{\infty}^{1,*-1}$ and so $(\Ker\Theta_M^G)^2=({}^I\!E_{\infty}^{1,*-1})^2=0$. So $\HH^*(A,M)^{G\uparrow}$ 
(which is isomorphic to $\HH^*(AG,M)$ if $M$ is an $AG$-algebra) is some singular extension of the algebra $\Im\Theta_M^G$ by the 
$\Im\Theta_M^G$-bimodule $\Ker\Theta_M^G$. On the other hand $\Im\Theta_M^G$ is a subalgebra of $\HH^*(A,M)^G$ and so we have an exact sequence of 
$\Im\Theta_M^G$-bimodules
$$
\Im\Theta_M^G\rightarrowtail\HH^*(A,M)^G\twoheadrightarrow \HH^*(A,M)^G/\Im\Theta_M^G.
$$
For simplicity we introduce the following notation: 
$$\AAA_M^G=\Im\Theta_M^G,\,X_M^G=\Ker\Theta_M^G,\,Y_M^G=\HH^*(A,M)^G/\Im\Theta_M^G,\,W_M^G=\Delta_{\rho}^{q-1}\HH^*(A,M)^H.$$
Note that $W_M^G\subset \AAA_M^G$. It is easy to check that if $C_q$ acts on some algebra $B$, then
\begin{equation}\label{derro}
\Delta_{\rho}(ab)=\Delta_{\rho}(a)b+a\Delta_{\rho}(b)-\Delta_{\rho}(a)\Delta_{\rho}(b)
\end{equation}
for any $a,b\in B$. Then it is easy to see that $W_M^G$ is an ideal of $\HH^*(A,M)^G$. Really, for $w\in\HH^*(A,M)^H$ and $f\in\HH^*(A,M)^G$ we obtain 
usinig the equality \eqref{derro} that $\Delta^{q-1}(w)f=\Delta^{q-1}(wf)$ and $f\Delta^{q-1}(w)=\Delta^{q-1}(fw)$. Let $\bar{\AAA}_M^G$ denote the 
quotient algebra $A_M^G/W_M^G$.

\section{Describing Hochschild cohomology}

During this section we assume that $A$ is a $\kk G$-module algebra and $M$ is an $(\mm,\hh)$-degenerated $A^eG$-module algebra.
Let $\bh$ denote the class of $\hh$ in $\HH^1(A,M)$. Then we have the following result.

\begin{lemma}\label{grcom}
The element $\bh\smile\bh$ lies in $\AAA_M^G$. Moreover, for any $\bar{a}\in(\AAA_M^G)_i$ the element $\bh\smile 
\bar{a}-(-1)^i\bar{a}\smile\bh$ lies in $\AAA_M^G$ too.
\end{lemma}
\begin{Proof} Let us construct ${}^I\!E$ as in Example \ref{cyclic}. We know that 
$${}^I\!E_{\infty}^{0,i}=(\AAA_M^G)_i={}^I\!E_{3}^{0,i}=\Ker{}^I\!d_2|_{{}^I\!E_2^{0,i}}.$$
Let us calculate ${}^I\!d_2(\bh\smile \bh)$. Using \eqref{derro} we obtain
$$\Delta_{\rho}(\hh\smile \hh)=-(\phi_{\mm}\smile \hh+\hh\smile \phi_{\mm}+\phi_{\mm}\smile \phi_{\mm})=-\delta_M^1(\mm\smile \hh-\hh\smile \mm+\mm\smile \phi_{\mm}).$$
 Direct calculations using formula \eqref{cyclic} show that
 $$
 \begin{aligned}
 \Delta_{\rho}^l(\mm\smile \hh)&=\Delta_{\rho}^l(\mm)\smile 
\hh+\sum\limits_{i=0}^l\sum\limits_{j=0}^{l-1}(-1)^{i+j-l}{{l}\choose{i}}{{i}\choose{l-j-1}}\Delta_{\rho}^i\mm\phi_{\Delta_{\rho}^j\mm},\\
  \Delta_{\rho}^l(\hh\smile 
\mm)&=\hh\smile\Delta_{\rho}^l\mm+\sum\limits_{i=0}^{l-1}\sum\limits_{j=0}^{l}(-1)^{i+j-l}{{l}\choose{i+1}}{{i+1}\choose{l-j}}\phi_{\Delta_{\rho}^i\mm}
\Delta_{\rho}^j\mm,\\
 \Delta_{\rho}^l(\mm\smile 
\phi_{\mm})&=\sum\limits_{i=0}^l\sum\limits_{j=0}^{l}(-1)^{i+j-l}{{l}\choose{i}}{{i}\choose{l-j}}\Delta_{\rho}^i\mm\phi_{\Delta_{\rho}^j\mm}.
   \end{aligned}
 $$
Since $\delta_M^0(\Delta_{\rho}^i\mm\Delta_{\rho}^j\mm)=\phi_{\Delta_{\rho}^i\mm}\Delta_{\rho}^j\mm+\Delta_{\rho}^i\mm\phi_{\Delta_{\rho}^j\mm}$, we have 
$\phi_{\Delta_{\rho}^i\mm}\Delta_{\rho}^j\mm=-\Delta_{\rho}^i\mm\phi_{\Delta_{\rho}^j\mm}$ in ${}^I\!E_2^{2,1}$. Note also that
 $$
 {{l}\choose{i}}{{i}\choose{l-j-1}}+{{l}\choose{i+1}}{{i+1}\choose{l-j}}+{{l}\choose{i}}{{i}\choose{l-j}}={{l+1}\choose{i+1}}{{i+1}\choose{l-j}}
 $$
 for all $i,j,l\in\mathbb{Z}$. So ${}^I\!d_2(\bh\smile \bh)$ is the class of
$$
\Delta_{\rho}^{q-1}(\mm\smile \hh-\hh\smile \mm+\mm\smile 
\phi_{\mm})=\sum\limits_{i=0}^{q-1}\sum\limits_{j=0}^{q-1}(-1)^{i+j+1-q}{{q}\choose{i+1}}{{i+1}\choose{q-1-j}}\Delta_{\rho}^i\mm\phi_{\Delta_{\rho}^j\mm}.
$$
in ${}^I\!E_2^{2,1}=H^2(C_q,\HH^1(A,M)^H)$. Since $p\mid{{q}\choose{i+1}}$ for $0\le i\le q-2$, we have
$$
{}^I\!d_2(\bh\smile \bh)=\sum\limits_{j=0}^{q-1}(-1)^{j}{{q}\choose{j+1}}\phi_{\Delta_{\rho}^j\mm}=0.
$$
So $\bh\smile \bh$ lies in $(\AAA_M^G)_2$.

Any $\bar{a}\in(\AAA_M^G)_i$ can be represented by such $a\in C^n(A,M)^H$ that $\delta_M^ia=0$ and $\Delta_{\rho}a=0$. Then 
$-\Delta_{\rho}(\hh\smile a)=\phi_{\mm}\smile a=\delta_M^i(\mm\smile a).$ So ${}^I\!d_2(\bh\smile \bar{a})$ can be represented by the element
$\Delta_{\rho}^{q-1}(\mm\smile a)=a$. Analogously ${}^I\!d_2(\bar{a}\smile\bh)$ can be represented by the element $(-1)^ia$. So $\bh\smile 
\bar{a}-(-1)^i\bar{a}\smile\bh$ lies in $\Ker{}^I\!d_2=\AAA_M^G$.
\end{Proof}

\begin{coro}\label{d2}
${}^I\!d_2\big((\bh\smile \bar{a})[0]\big)=\bar{a}[2]$ for any $\bar{a}\in\AAA_M^G$.
\end{coro}
\begin{Proof}
The statement follows from the proof of Lemma \ref{grcom}.
\end{Proof}

Let $A$ be a $\kk G$-module algebra and $M$ be an $(\mm,\hh)$-degenerated $A^eG$-module algebra.
Let $D_{\hh}:\AAA_M^G\rightarrow\AAA_M^G$ be a map defined by the equality $D_{\hh}(a)=\bh\smile a-(-1)^ia\smile \bh$ for $a\in(\AAA_M^G)_i$. It is easy to see that 
$D_{\hh}$ sends $W_M^G$ to $W_M^G$. We denote the induced map from $\bar{\AAA}_M^G$ to itself by $D_{\hh}$ too. Note that $D_{\hh}$ is a graded derivation of 
$\bar{\AAA}_M^G$, i.e. $D_{\hh}(a\smile b)=D_{\hh}(a)\smile b+(-1)^ia\smile D_{\hh}(b)$ for $a\in(\bar{\AAA}_M^G)_i$, $b\in\bar{\AAA}_M^G$.

\begin{Def}
Let $B$ be a graded algebra and $D:B\rightarrow B$ be a graded derivation. Let us define the graded algebra $B[x,D]$. Its underlining graded vector 
space is the space $B[x]$ with $n$-th component generated by the elements of the form $bx^i$ where $b\in B_{n-i}$. The multiplication of two elements 
from $B\subset B[x,D]$ or two elements from $\kk[x]\subset B[x,D]$ is defined as usual. We define the left multiplication of $b\in B_i$ by the elemet 
$x$ by the equality $xb=(-1)^ibx+D(b)$.
\end{Def}

Now we are ready to almost describe the algebra $\HH^*(A,M)^G$.

\begin{theorem}\label{Gstabdes}
Let $A$ be a $\kk G$-module algebra and $M$ be an $(\mm,\hh)$-degenerated $A^eG$-module algebra.
Then there is an isomorphism of graded algebras $$\HH^*(A,M)^G/W_M^G\cong\bar{\AAA}_M^G[x,D_{\hh}]/\langle x^2-\bh\smile\bh\rangle.$$
\end{theorem}

Let $B$ be a graded algebra and $L$ be a graded $B$-bimodule. Then we write $L[t]$ for the graded $B$-bimodule, which equals to $L$ as vector space, 
has the grading defined by the equality $L[t]_n=L_{n+t}$ and has the bimodule structure defined by the equality $al[t]b=(-1)^{it}(alb)[t]$ for $a\in B$, 
$l\in L$ and $b\in B_i$, where $l[t]$ denotes the element of $L[t]$ corresponding to $l$. To prove Theorem \ref{Gstabdes} we need some lemmas.

\begin{lemma}\label{hmul}
Let $\psi(\bar{a})$ ($\bar{a}\in\AAA_M^G$) be the class of the element $\bar{a}\smile\bh\in\HH^*(A,M)^G$ in $Y_M^G$.
Then $\psi:\AAA_M^G\rightarrow Y_M^G[1]$ is a homorphism of graded $\AAA_M^G$-bimodules. Moreover, $\Ker\psi=W_M^G$.
\end{lemma}
\begin{Proof} By Lemma \ref{grcom} for any $\bar{a}\in\AAA_M^G$ and $\bar{b}\in(\AAA_M^G)_i$ there is such $\bar{c}\in\AAA_M^G$ that
$$(\bar{a}\smile\bar{b})\smile\bh=\bar{a}\smile(\bar{b}\smile\bh)=\bar{a}\smile((-1)^i\bh\smile\bar{b}+\bar{c})=(-1)^i(\bar{a}\smile 
\bh)\smile\bar{b}+\bar{a}\bar{c}.$$
So we have $\psi(\bar{a}\smile\bar{b})=\bar{a}\smile\psi(\bar{b})=\psi(\bar{a})\smile\bar{b}$, i.e. $\psi$ is an $\AAA_M^G$-bimodule homomorphism.

Let now take some $\bar{a}\in(\AAA_M^G)_i$. The class of $\bar{a}\smile\bh$ equals to zero in $Y_M^G$ iff $\bar{a}\smile\bh\in A_M^G$. But 
$\bar{a}\smile\bh\in A_M^G$ is equivalent to $0={}^I\!d_2\big((\bar{a}\smile\bh)[0]\big)=(-1)^i\bar{a}[2]$ (see Corollary \ref{d2}). But the class of 
$\bar{a}[2]$ equals to zero in ${}^I\!E_2^{2,i}=\Ho^2(C_q,\HH^i(A,M)^H)$ iff  there is such $\bar{b}\in\HH^i(A,M)^H$ that 
$\Delta_{\rho}^{q-1}\bar{b}=\bar{a}$. But the last condition is eqivalent to the fact that $\bar{a}$ lies in $W_M^G$. We obtain that $\Ker\psi=W_M^G$.
\end{Proof}

The next lemma says that the homomorphism $\psi$ from the previous lemma is surjective. This fact can be deduced from dimension conditions, but we give 
a proof that do not use arguments related to finiteness of any dimension.

\begin{lemma}\label{sur}
Any element of $\HH^*(A,M)^G$ can be represented in the form $\bar{a}+\bar{b}\smile\bh$, where $\bar{a},\bar{b}\in\AAA_M^G$.
\end{lemma}
\begin{Proof}
We prove using induction that $\HH^i(A,M)^G=(\AAA_M^G)_i+(\AAA_M^G)_{i-1}\smile\bh$. This statemen is true for $i=0$, since 
$\HH^0(A,M)^G=(\AAA_M^G)_0$. Suppose that we know that our statement is true for $i=n$. Let us take some $f\in\HH^{n+1}(A,M)^G$.
We have ${}^I\!d_2(f[0])=g[2]$ for some $g\in\HH^n(A,M)^G$. By induction hypothesis we obtain $g=\bar{a}+\bar{b}\smile\bh$. Since 
${}^I\!d_2(g[2])=0$, we have ${}^I\!d_2(g[0])=0$ and so $\bar{b}\in W_M^G$ by Lemma \ref{hmul}. Then $g\in(\AAA_M^G)_n$ and we have 
${}^I\!d_2(f[0]-(-1)^n(g\smile\bh)[0])=0$, i.e. $f-(-1)^ng\smile\bh$ lies in $(\AAA_M^G)_{n+1}$. The lemma is proved.
\end{Proof}

\begin{Proof}[Proof of Theorem \ref{Gstabdes}.] Let us define the map $$\Phi:\bar{\AAA}_M^G[x,D_{\hh}]/\langle x^2-\bh\smile \bh\rangle\rightarrow 
\HH^*(A,M)^G/W_M^G$$ by the formula $\Phi\big((\bar{a}+W_M^G)+(\bar{b}+W_M^G)x\big)=\bar{a}+\bar{b}\smile\bh+W_M^G$ for 
$\bar{a},\bar{b}\in\AAA_M^G$. It is clear that the definition of $\Phi$ is correct and it is not hard to check that $\Phi$ is a homomorphism of 
algebras. It follows from Lemma \ref{sur} that $\Phi$ is surjective. Assume now that 
$$0=\Phi\big((\bar{a}+W_M^G)+(\bar{b}+W_M^G)x\big)=\bar{a}+\bar{b}\smile\bh+W_M^G,$$
 i.e. $\bar{a}+\bar{b}\smile\bh\in W_M^G$. Then $\bar{b}\smile\bh\in\AAA_M^G$ and so $\bar{b}\in W_M^G$ by Lemma \ref{Gstabdes}. But then we 
obtain that $\bar{a}\in W_M^G$ too and so $(\bar{a}+W_M^G)+(\bar{b}+W_M^G)x=0$. We prove the injectivity of $\Phi$ and so $\Phi$ is an isomorphism of 
algebras.
\end{Proof}

\begin{rema}
If $M=A$ or $M=AG$, then the algebra $\HH^*(A,M)^G$ is graded commutative (see Section \ref{GActH}). In this case we have $\HH^*(A,M)^G/W_M^G\cong\bar{\AAA}_M^G[x,0]/\langle 
x^2-\bh\smile\bh\rangle$. If additionally we have $\bh\smile\bh=0$ (this equality is automatically true if $\HH^*(A,M)$ is graded 
commutative and $\charr\kk\not=2$), then $\HH^*(A,M)^G/W_M^G\cong\bar{\AAA}_M^G[x,0]/\langle x^2\rangle$ is a trivial singular extension of 
$\bar{\AAA}_M^G$ by $\bar{\AAA}_M^G[-1]$.
\end{rema}

\section{Comparing a kernel and a cokernel of $\Theta_M^G$}

In this section we assume again that $A$ is a $\kk G$-module algebra and $M$ is an $(\mm,\hh)$-degenerated $A^eG$-module algebra.
The aim of this section is to prove the following result, which  gives some additional relation between $Y_M^G$ and $X_M^G$ (and so between 
$\HH^*(A,M)^G$ and $\HH^*(A,M)^{G\uparrow}$).

\begin{theorem}\label{filtering}
Let $A$ be a $\kk G$-module algebra and $M$ be an $(\mm,\hh)$-degenerated $A^eG$-module algebra.
Then there are two filtrations of graded $\AAA_M^G$-bimodules $$0=X_{q-1}\subset X_{q-2}\subset\cdots\subset X_1\subset X_0=X_M^G$$ and 
$$0=Y_{q-1}\subset Y_{q-2}\subset\cdots\subset Y_1\subset Y_0=Y_M^G$$ such that $Y_{i}/Y_{i+1}\cong X_{q-2-i}/X_{q-1-i}$ for $0\le i\le q-2$. Moreover, 
 if $\Delta_{\rho}^i\HH^*(A,M)^H=\HH^*(A,M)^G$ for some $0\le i\le q-1$, then $X_{q-1-i}=0=Y_{i+1}$ and $X_{q-2-i}=X\cong Y=Y_i$.
\end{theorem}

On one hand, Theorem \ref{filtering} says that there is some relation between $X_M^G$ and $Y_M^G$, but, on the other hand, it says that this relation 
is very complicated in general.
To prove this result we need some auxilary facts.

\begin{lemma}\label{cycnull}
Let $f\in C^j(A,M)^H$ be such that $\delta_M^jf=0$ and $\Delta_{\rho}^if=0$ for some $0\le i\le q$. Then there are $f_0\in C^{j}(A,M)^H$ and $f_1\in 
C^{j+1}(A,M)^H$ such that
$\Delta_{\rho}^{q-i}f_0=f$, $\delta_M^jf_0=\Delta_{\rho}^if_1$ and $\delta_M^{j+1}f_1=0$.
\end{lemma}
\begin{Proof} If $i=0$ or $i=q$, then the lemma is obvious. Assume that $0<i<q$. Let us take $f_1=(-1)^{j-1}f\smile \hh$ and
$$
f_0=\sum\limits_{k=0}^{i-1}\sum\limits_{l=0}^{i-1}(-1)^{l+k-i-1}{{i}\choose{k}}{{k}\choose{i-l-1}}\Delta_{\rho}^kf\Delta_{\rho}^l\mm.
$$
It is clear that $\delta_M^{j+1}f_1=0$. Easy computation using \eqref{derro} shows that $\delta_M^jf_0=\Delta_{\rho}^if_1$. The equality 
$\Delta_{\rho}^{q-i}f_0=f$ can be obtained by direct computation too. Since it is more complicated, let us perform it. Using \eqref{derro} we obtain 
that $\Delta_{\rho}^{q-i}f_0$ equals to
\begin{multline*}
\sum\limits_{k=0}^{i-1}\sum\limits_{l=0}^{i-1}(-1)^{l+k-i-1}{{i}\choose{k}}{{k}\choose{i-l-1}}\sum\limits_{s=0}^{q-i}\sum\limits_{t=0}^{q-i}(-1)^{s+t-q+i
}{{q-i}\choose{s}}{{s}\choose{q-i-t}}\Delta_{\rho}^{k+s}f\Delta_{\rho}^{l+t}\mm\\
=\sum\limits_{s=0}^{i-1}\sum\limits_{t=0}^{q-1}(-1)^{s+t+q-1}\left(\sum\limits_{k=0}^s\sum\limits_{l=0}^t{{i}\choose{k}}{{k}\choose{i-l-1}}{{q-i}
\choose{s-k}}{{s-k}\choose{q-i-t+l}}\right)\Delta_{\rho}^sf\Delta_{\rho}^t\mm.
\end{multline*}
For any integers $n$, $s$ and $l$ such that $0\le l\le n$ we have
$
\sum\limits_{k=0}^s{{l}\choose{k}}{{n-l}\choose{s-k}}={{n}\choose{s}}.
$ To prove this, it is enough to calculate the coefficient of $x^s$ on both sides of the eqaulity $(x+1)^n=(x+1)^l(x+1)^{n-l}$. Then for $0\le s\le i-1$ and $0\le t\le q-1$ we have
\begin{multline*}
\sum\limits_{k=0}^s\sum\limits_{l=0}^t{{i}\choose{k}}{{k}\choose{i-l-1}}{{q-i}\choose{s-k}}{{s-k}\choose{q-i-t+l}}\\
=\sum\limits_{k=0}^s{{i}\choose{k}}{{q-i}\choose{s-k}}\sum\limits_{l=0}^t\left({{k}\choose{i-l-1}}{{s-k}\choose{q-i-t+l}}\right)\\
=\sum\limits_{k=0}^s{{i}\choose{k}}{{q-i}\choose{s-k}}{{s}\choose{q-t-1}}={{q}\choose{s}}{{s}\choose{q-t-1}}.
\end{multline*}
Since $p\mid {{q}\choose{s}}$ for $0<s<q$, we have
$
\Delta_{\rho}^{q-i}f_0=f.
$
\end{Proof}

\begin{coro}\label{bounull}
Let $f\in C(A,M)^H$ be such that $\delta_M\Delta_{\rho}^if=0$ for some $0\le i\le q$. Then there is $f_0\in C(A,M)^H$ such that 
$\Delta_{\rho}^{q-i}f_0=\delta_Mf$ and $\delta_Mf_0=0$.
\end{coro}
\begin{Proof}
It follows from Lemma \ref{cycnull} that there are such $g_0,g_1\in C(A,M)^H$ that $\Delta_{\rho}^{i}g_0=\Delta_{\rho}^{i}f$, 
$\delta_Mg_0=\Delta_{\rho}^{q-i}g_1$ and $\delta_Mg_1=0$. Since $\Delta_{\rho}^i(f-g_0)=0$ and $C(A,M)^H$ is a projective $C_q$-module, we have some 
$g_2\in C(A,M)^H$ such that $\Delta_{\rho}^{q-i}g_2=f-g_0$. Then we can take $f_0=g_1+\delta_Mg_2$.
\end{Proof}

\begin{Proof}[Proof of Theorem \ref{filtering}.]
Let $\pi$ denote the canonical projection $\HH^*(A,M)^G\rightarrow Y_M^G$.
Let us define $Z_i$ ($0\le i\le q-1$) by the equality $Z_i=\AAA_M^G\cap\Delta_{\rho}^i\HH^*(A,M)^H$. Let introduce $Y_i=\pi(\bh\smile Z_i)$ for 
$0\le i\le q-1$. Since $Z_{q-1}=W_M^G$ and $Z_{0}=\AAA_M^G$, it follows from Lemmas \ref{hmul} and \ref{sur} that $Y_{q-1}=0$ and $Y_0=Y_M^G$.
It is clear that $Z_i$ is a graded $\AAA_M^G$-subbimodule of $\AAA_M^G$ and so $Y_i$ is a graded $\AAA_M^G$-subbimodule of $Y_M^G$ for any $0\le i\le 
q-1$. Moreover, if $\Delta_{\rho}^i\HH^*(A,M)^H=\HH^*(A,M)^G$ for some $i$, then $Z_{i+1}=0$, $Z_i=\AAA_M^G$ and so $Y_{i+1}=0$ and $Y_i=Y_M^G$.

Let now define $X_i$ ($0\le i\le q-1$) as a set of elements $\bar{a}\in X_M^G$ that can be represented by an element from $C(A,M)^G\cap 
\delta_M\Delta_{\rho}^iC(A,M)^H$. By definition of $X_M^G$ we have $X_0=X_M^G$ and $X_{q-1}=0$.
Moreover, $X_i$ is a graded $\AAA_M^G$-subbimodule of $X_M^G$ since $b\smile (\delta_M\Delta_{\rho}^if)\smile c=(-1)^k\delta_M\Delta_{\rho}^i(b\smile 
f\smile c)$ for $b\in C^k(A,M)^G$ and $c\in C^l(A,M)^G$ such that $\delta_M^kb=\delta_M^lc=0$.

It remains to prove that $Y_{i}/Y_{i+1}\cong X_{q-2-i}/X_{q-1-i}$ for $0\le i\le q-2$. Let us take some $\bar{y}\in Y_i$. Then we have 
$\bar{y}=\pi(\bh\smile\bar{a})$, where $\bar{a}$ is the class of some element $a\in C(A,M)^G$ of the form $a=\Delta_{\rho}^if+\delta_Mg$, where 
$f,g\in C(A,M)^H$ and $\delta_Mf=0$. Then we have $\delta_M\Delta_{\rho}^{q-i}g=0$. By Corollary \ref{bounull} there is some $g_0\in C(A,M)^H$ 
satisfying the conditions $\delta_Mg_0=0$ and $\Delta_{\rho}^ig_0=\delta_Mg$. Let us consider an element $u=a\smile\Delta_{\rho}^{q-1-i}\mm-f-g_0$. It is 
easy to see that $\Delta_{\rho}^iu=0$ and so there is such $v\in C(A,M)^H$ that $\Delta_{\rho}^{q-i}v=u$.
It is easy to see that the element $\delta_M\Delta_{\rho}^{q-2-i}(a\smile \mm-\Delta_{\rho}v)$ lies in $C(A,M)^G\cap 
\delta_M\Delta_{\rho}^{q-2-i}C(A,M)^H$. Let $\bar{v}$ denote the class of $\delta_M\Delta_{\rho}^{q-2-i}(a\smile \mm-\Delta_{\rho}v)$ in $X_{q-2-i}$.
We define the map $\Phi_i:Y_{i}/Y_{i+1}\rightarrow X_{q-2-i}/X_{q-1-i}$ by the formula
$$
\Phi_i(\bar{y}+Y_{i+1})=\bar{v}+X_{q-1-i}.
$$

Let us prove the correctness of the definition of $\Phi_i$. It is enough to prove that $\bar{v}\in X_{q-1-i}$  for any choice of intermediate elements 
if $\bar{y}\in Y_{i+1}$. Let $\bar{a}$, $a$, $f$, $g$, $g_0$, $u$ and $v$ be the elemets used in construction of 
$\bar{v}$ above. Since $\bar{y}\in Y_{i+1}$, we have $y=\pi(\bh\smile\bar{a}')$, where the elemet $\bar{a}'$ can be represented by $a'\in C(A,M)^G$ 
of the form $a'=\Delta_{\rho}^{i+1}f'+\delta_Mg'$, where $f',g'\in C(A,M)^H$ and $\delta_Mf'=0$. We have $\bar{a}-\bar{a}'\in W_M^G$ by Lemma 
\ref{hmul} and so 
$$a=a'+\Delta_{\rho}^{q-1}x+\delta_My=\Delta_{\rho}^{i+1}(f'+\Delta_{\rho}^{q-2-i}x)+\delta_M(g'+y),$$
 where $x,y\in C(A,M)^H$ and $\delta_Mx=0$. Let $w=f'+\Delta^{q-2-i}x$, $z=g'+y$. Then we have 
$a=\Delta_{\rho}^if+\delta_Mg=\Delta_{\rho}^{i+1}w+\delta_Mz$. Since $\delta_M\Delta_{\rho}^{q-1-i}z=0$, Corollary \ref{bounull} gives us an element 
$z_0$ such that $\delta_Mz_0=0$ and $\Delta_{\rho}^{i+1}z_0=\delta_Mz$. Let consider the element $r=f+g_0-\Delta_{\rho}(w+z_0)$. Since 
$\Delta_{\rho}^ir=0$, there is $r_0\in C(A,M)^H$ such that $\Delta_{\rho}^{q-i}r_0=r$. Then 
 $$\bar{v}_0=\delta_M\Delta_{\rho}^{q-2-i}\big(a\smile \mm-\Delta_{\rho}(v+r_0)\big)$$
 represents some element of $X_{q-2-i}$ and, moreover, $\bar{v}+X_{q-1-i}=\bar{v}_0+X_{q-1-i}$. Since
 $$
\delta_M\Delta_{\rho}^{q-2-i}\big(a\smile \mm-\Delta_{\rho}(v+r_0)\big)=\delta_M\big(a\smile 
\Delta_{\rho}^{q-2-i}\mm-\Delta_{\rho}^{q-1-i}(v+r_0)-w-z_0\big)
 $$
and $\Delta_{\rho}\big(a\smile \Delta_{\rho}^{q-2-i}\mm-\Delta_{\rho}^{q-1-i}(v+r_0)-w-z_0\big)=0$, we have $\bar{v}_0=0$ in $X_M^G$. So $\bar{v}\in 
X_{q-1-i}$ and the correctness is proved.

Let $a\in C^k(A,M)^G$, $b\in C^l(A,M)^G$ be such that $\delta_M^ka=\delta_M^lb=0$. By Lemma \ref{grcom} there are $c\in C^{l+1}(A,M)^G$ and $u\in 
C^{l}(A,M)^H$ such that $\delta_M^lc=0$ and
$$\hh\smile b=(-1)^lb\smile \hh+c+\delta_M^{l}u.$$
Then we have
\begin{multline*}
\delta_M^{k+l}\Delta_{\rho}^{q-2-i}(a\smile \mm\smile b)=(-1)^ka\smile\phi_{\Delta_{\rho}^{q-2-i}\mm}\smile b=(-1)^{k+1}\Delta_{\rho}^{q-1-i}(a\smile \hh\smile 
b)\\
=(-1)^{k+l+1}\Delta_{\rho}^{q-1-i}(a\smile b\smile \hh)-(-1)^k\Delta_{\rho}^{q-1-i}\big(a\smile (c+\delta_M^{l}u)\big)\\
=\delta_M^{k+l}\Delta_{\rho}^{q-2-i}(a\smile b\smile \mm)-\delta_M^{k+l}\Delta_{\rho}^{q-1-i}(a\smile u).
\end{multline*}
for $a\in C^k(A,M)^G$, $b\in C^l(A,M)^G$  such that $\delta_M^ka=\delta_M^lb=0$.
Now it is easy to see that $\Phi_i$ is a homomorphism of $\AAA_M^G$-modules.

Let us prove that $\Phi_i$ is injective. Assume that $\bar{y}\in Y_i$ is such that $\Phi_i(\bar{y})=0$. Let $\bar{a}$, $a$, $f$, $g$, $g_0$, $u$, $v$ and 
$\bar{v}$ be as above. Then $\Phi_i(\bar{y})=0$ iff $\bar{v}\in X_{q-1-i}$. This means that
$$
\delta_M\Delta_{\rho}^{q-2-i}(a\smile \mm-\Delta_{\rho}v)=\delta_M\Delta_{\rho}^{q-1-i}w+\delta_Mz
$$
for some $w\in C(A,M)^H$ and $z\in C(A,M)^G$. Then $a\smile {\Delta_{\rho}^{q-2-i}\mm}-\Delta_{\rho}^{q-1-i}(v+w)-z$ represents some element of 
$\HH^*(A,M)$. On the othe hand,
$$
\Delta_{\rho}^{i+1}\big(a\smile {\Delta_{\rho}^{q-2-i}\mm}-\Delta_{\rho}^{q-1-i}(v+w)-z\big)=a.
$$
So $\bh\smile\bar{a}\in Z_{i+1}$, i.e. $\bar{y}$ represents zero element of $Y_i/Y_{i+1}$.

It remains to prove that $\Phi_i$ is surjective. Let us take some $\bar{x}\in X_{q-2+i}$. By definition $\bar{x}$ can be represented by 
$\delta_M\Delta^{q-2-i}a$ for some $a\in C(A,M)^H$ such that $\delta_M\Delta^{q-1-i}a=0$.
Let consider $\bar{y}=\pi(\bh\smile \bar{a})$, where $\bar{a}$ is the class of $\Delta_{\rho}^{q-1}a$ in $\AAA_M^G$. Since 
$\Delta_{\rho}^{q-1}a=\Delta_{\rho}^i(\Delta_{\rho}^{q-1-i}a)$, it is clear that $\bar{y}\in Y_i$ and
$\Phi_i(\bar{y}+Y_{i+1})=\bar{v}+X_{q-1-i}$, where $\bar{v}$ is the class of $\delta_M\Delta_{\rho}^{q-2-i}(\Delta_{\rho}^{q-1}a\smile 
\mm-\Delta_{\rho}v)$ in $X_{q-2-i}$ for some $v\in C(A,M)^H$ such that 
$$\Delta_{\rho}^{q-i}v=\Delta_{\rho}^{q-1}a\smile \Delta_{\rho}^{q-1-i}\mm-\Delta_{\rho}^{q-1-i}a$$
(see the construction above). But then $\bar{x}-\bar{v}$ is the class of
$$
\delta_M\Delta^{q-2-i}(a-\Delta_{\rho}^{q-1}a\smile \mm+\Delta_{\rho}v)
$$
in $X_{q-2-i}$. Since
$$\Delta_{\rho}\Delta^{q-2-i}(a-\Delta_{\rho}^{q-1}a\smile \mm+\Delta_{\rho}v)=\Delta_{\rho}^{q-1-i}a-\Delta_{\rho}^{q-1}a\smile 
\Delta_{\rho}^{q-1-i}\mm+\Delta_{\rho}^{q-i}v=0,$$
we have $\bar{x}=\bar{v}$ and so $\bar{x}+X_{q-1-i}=\Phi_i(\bar{y})$.
\end{Proof}

\section{The case of a Hopf-Galois covering}

In this section we consider the case where $A$ is a $\kk G$-covering of some algebra $B$, i.e. there is some $G$-graded algebra $B$ such that $A=BG^*$ 
(see Example \ref{covering}). The $G$-action on such $A$ is defined in Section \ref{GActH}. It is clear that $A\cong(\kk G)^{\dim_{\kk}B}$ as $\kk 
G$-module. In particular, we have $1_A=\Delta_{\rho}^{q-1}\mm$, where $\mm$ denotes the element $\sum\limits_{\alpha\in H}p_{\alpha}\in A^H$. On the other 
hand it is possible that there is no such $\hh\in\Der(A,A)^H$ that $A$ is an $(\mm,\hh)$-degenerated algebra. It is shown in the following example.

\begin{Ex}
Let $\charr\kk=2$, $G=C_2$ and $B=\kk[x,y]/\langle x^2-y^2, xy\rangle$. Let introduce such $G$-grading on $B$ that $B_{\rho}={}_{\kk}\langle x\rangle$. 
Then $A=BG^*\cong \kk Q/I$, where $(Q,I)$ is the quiver
$Q=\xy
  (25,10)*+\txt{\scriptsize 0};
  (25,10)*+\txt{\scriptsize 0}
  **\crv{(30, 5)&(25,0)&(20, 5)}
  ?>*\dir{>}
  ?(.5) *!LD!/^-6pt/\txt{\tiny $y_0$};
    (25,10)*+\txt{\scriptsize 0};
    (45,10)*+\txt{\scriptsize 1}
**\crv{(35, 15)}
  ?>*\dir{>}
  ?(.5) *!LD!/^-2pt/\txt{\tiny $x_0$};
      (45,10)*+\txt{\scriptsize 1};
    (25,10)*+\txt{\scriptsize 0}
**\crv{(35, 5)}
  ?>*\dir{>}
  ?(.5) *!LD!/^-8pt/\txt{\tiny $x_1$};
(45,10)*+\txt{\scriptsize 1};
(45,10)*+\txt{\scriptsize 1}
  **\crv{(50, 5)&(45,0)&(40, 5)}
  ?>*\dir{>}
  ?(.5) *!LD!/^-6pt/\txt{\tiny $y_1$};
  \endxy$ with ideal of relations $I=\langle x_1x_0-y_0^2,x_0x_1-y_1^2,x_0y_0,y_1x_0,x_1y_1,y_0x_1\rangle$. The $G$-action on $A$ is defined by the 
equalities $$\rho(x_i)=x_{i+1},\,\rho(y_i)=y_{i+1}, \,\rho(e_i)=e_{i+1}\,(i\in\mathbb{Z}_2),$$
  where $e_i$ is an idempotent corresponding to a vertex $i$ of the quiver $Q$. Then we can take $\mm=e_0$ and it is easy to check that 
$\Delta_{\rho}\hh\not=\phi_{e_0}$ for any $\hh\in\Der(A,A)^H$.
\end{Ex}

Nonetheless, there is a case in which the algebra $A=BG^*$ is $(\mm,\hh)$-degenerated. We say that a $G$-grading on $B$ is $C_q$-singular if 
$B_{\alpha}B_{\beta}=0$ for $\alpha,\beta\in G$, $\alpha,\beta\not\in H$.

\begin{lemma}
Assume that the $G$-grading on $B$ is $C_q$-singular. Then the algebra $A=BG^*$ is an $(\mm,\hh)$-degenerated $A^eG$-module algebra. Moreover, the map 
$\hh\in\Der(A,A)^H$ can be choosen in such a way that $\hh\smile \hh=0$ in $C^2(A,A)$.
\end{lemma}
\begin{Proof} Let us take $\mm=\sum\limits_{\alpha\in H}p_{\alpha}$ and
$$\hh(ap_{\rho^j\beta})=\begin{cases}
ap_{\rho^j\beta},&\mbox{if $i<j$,}\\
0,&\mbox{otherwise}\end{cases}\,(a\in B_{\rho^i\alpha},0\le i,j\le q-1,\alpha,\beta\in H).$$ 
Direct calculations show that $\hh\in\Der(A,A)^H$ and the equalities $(\Delta_{\rho}\hh)(ap_{\alpha})=\mm ap_{\alpha}-ap_{\alpha}\mm$ and 
$\hh(ap_{\alpha})\hh(bp_{\beta})=0$ are true for all $a,b\in B$ and $\alpha,\beta\in G$.
\end{Proof}

Since the algebras $\HH^*(B)$, $\HH^*(A)$ and $\HH^*(A,AG)^G$ are graded commutative we have the following corollaries.

\begin{coro}
Let $B$ be equipped with a $C_q$-singular $G$-grading.
Then there are a graded subalgebra $\AAA_{AG}^G$ of $\HH^*(A,AG)^G$ and an ideal $X_{AG}^G$ of $\HH^*(B)$ such that
$$
(X_{AG}^G)^2=0\mbox{ and }\HH^*(B)/X_{AG}^G\cong\AAA_{AG}^G.
$$
Moreover, there is also an ideal $W_{AG}^G\subset\AAA_{AG}^G$ of $\HH^*(A,AG)^G$ such that
$$
\HH^*(A,AG)^G/W_{AG}^G\cong (\AAA_{AG}^G/W_{AG}^G)[x,0]/x^2
$$
and for any $i\ge 0$ the equality $\dim_{\kk}(W_{AG}^G)_i+\dim_{\kk}(X_{AG}^G)_{i+1}=\dim_{\kk}(\AAA_{AG}^G)_i$ holds.
\end{coro}

\begin{coro}
Let $B$ be equipped with a $C_q$-singular $G$-grading.
Then there are a graded subalgebra $\AAA_A^G$ of $\HH^*(A)^G$ and an ideal $X_{A}^G$ of $\HH^*(B)_0$ such that
$$
(X_{A}^G)^2=0\mbox{ and }\HH^*(B)_0/X_{A}^G\cong\AAA_{A}^G.
$$
Moreover, there is also an ideal  $W_{A}^G\subset\AAA_{A}^G$ of $\HH^*(A)^G$ such that
$$
\HH^*(A)^G/W_{A}^G\cong (\AAA_{A}^G/W_{A}^G)[x,0]/x^2
$$
and for any $i\ge 0$ the equality $\dim_{\kk}(W_{A}^G)_i+\dim_{\kk}(X_{A}^G)_{i+1}=\dim_{\kk}(\AAA_{A}^G)_i$ holds.
\end{coro}

Let $\RR$ be a finite dimensional $\kk$-algebra. We denote by $D\RR$ the $\RR$-bimodule $\Hom_{\kk}(\RR,\kk)$. Let $T\RR$ be a trivial singular extension of $\RR$ 
by $D\RR$, i.e. $T\RR$ is an algebra, whose underlining space is $\RR\oplus D\RR$ and the multiplication is defined by the equality
$(a,\hat{u})(b,\hat{v})=(ab,a\hat{v}+\hat{u}b)$ for $a,b\in \RR$ and $\hat{u},\hat{v}\in D\RR$. There is a $\mathbb{Z}$-grading on $T\RR$ such that $(T\RR)_0=\RR$, 
$(T\RR)_1=D\RR$ and $(T\RR)_i=0$ for $i\not\in\{0,1\}$. The algebra $\hat \RR=T\RR\#\kk\mathbb{Z}^*$ is called the repetitive algebra of $\RR$. Since the notion of 
$\kk G^*$-algebra was introduce only for finite group $G$, let us explain what does $T\RR\#\kk\mathbb{Z}^*$ means. As a space it is simply 
$\oplus_{i\in\mathbb{Z}}T\RR p_i$. The multiplication is defined by the equality $xp_iyp_j=xy_{i-j}p_j$ ($x,y\in T\RR$, $i,j\in\mathbb{Z}$).

 On the other hand, for any $n\ge 1$ there is a group epimorphism $\pi_n:\mathbb{Z}\rightarrow C_n$ and it defines a $C_n$-grading on $T\RR$. Let $\RR_n$ 
denote the algebra $T\RR\#\kk C_n^*$. In particular, we have $\RR_1=T\RR$. There is a linear map $\ve:T\RR\rightarrow\kk$ that sends $(a,\hat{u})$ ($a\in \RR$, 
$\hat{u}\in D\RR$) to $\hat{u}(1)$. This map defines a nondegenerated symmetric associative bilinear form $\langle,\rangle:T\RR\times T\RR\rightarrow\kk$ by 
the equality $\langle x, y\rangle=\ve(xy)$ ($x,y\in T\RR$). So $T\RR$ is a symmetric algebra.
On the other hand $\ve$ induces linear maps $\ve_n:\RR_n\rightarrow \kk$ ($n\ge 1$) and $\ve_{\infty}:\hat{\RR}\rightarrow\kk$ by the rules 
$\ve_n(xp_{\alpha})=\ve(x)$ and $\ve_{\infty}(xp_i)=\ve(x)$ for $x\in T\RR$, $\alpha\in C_n$ and $i\in\mathbb{Z}$. These maps define nondegenerated 
associative bilinear forms $\langle,\rangle_n:\RR_n\times \RR_n\rightarrow\kk$ ($n\ge 1$) and $\langle,\rangle_{\infty}:\hat{\RR}\times 
\hat{\RR}\rightarrow\kk$. So $\RR_n$ ($n\ge 1$) and $\hat{\RR}$ are Frobenius algebras. They have Nakayama automorphisms $\nu_n\in\Aut_{\kk}(\RR_n)$ ($n\ge 1$) and 
$\nu_{\infty}\in\Aut_{\kk}(\hat{\RR})$ defined by the equalities $\nu_n(xp_{\alpha})=xp_{\pi_n(-1)\alpha}$ and $\nu_{\infty}(xp_i)=xp_{i-1}$ for $x\in T\RR$, 
$\alpha\in C_n$ and $i\in\mathbb{Z}$. In some sense we can interpret $\RR_n$ as $\hat{\RR}/\langle\nu_{\infty}^n\rangle$.

Since the action of $C_n$ on $\RR_n$ is defined by a map $\eta:C_n\rightarrow \Aut_{\kk}(\RR_n)$ that sends some generator of $C_n$ to $\nu_n$, it follows from 
\cite[Corollary 2]{Volk6} that 
\begin{equation}\label{nustab}
\HH^*(\RR_n)=\HH^*(\RR_n)^{C_n}.
\end{equation}
In particular, if $\charr\kk\nmid n$, then $\HH^*(\RR_n)\cong\HH^*(T\RR)_0$.

If $\charr\kk=p$, $q$ is a positive power of $p$ and $n=qs$, where $p\nmid s$, then we have an exact sequence of groups $C_s\rightarrowtail C_n\twoheadrightarrow 
C_q$ and so the results of this paper can be applied to this case. In view of \eqref{nustab} our formulas are simplified and we obtain the following 
theorem. 

\begin{theorem}\label{TR}
Let $\RR$ be a finite dimensional $\kk$-algebra and $n\ge 1$ be some integer such that $\charr \kk\mid n$. Then there is a graded algebra $\AAA_n$ such 
that $\HH^*(\RR_n)$ is a trivial singular extension of $\AAA_n$ by $\AAA_n[-1]$ and $\HH^*(T\RR)_0$ is some singular extension of $\AAA_n$ by $\AAA_n[-1]$. 
In other words,
$$
\HH^*(\RR_n)\cong \AAA_n[x,0]/x^2
$$
and there is an exact sequence of $\AAA_n$-modules
$$
\AAA_n[-1]\rightarrowtail \HH^*(T\RR)_0\twoheadrightarrow \AAA_n,
$$
which becomes an exact sequence of algebras if we equip $\AAA_n[-1]$ with zero multiplication.
\end{theorem}

\noindent{{\bf Addresses:}
\newline 
Eduardo do Nascimento Marcos\\
IME-USP(Dept de Matem\'atica )
\newline
Rua de Mat\~ao 1010, Cidade Universit\'aria, S\~ao Paulo-SP
\newline
e-mail:  enmarcos@ime.usp.br
\newline
Yury Volkov \\
Temporary:
IME-USP(Dept de Matem\'atica )
\newline
Rua de Mat\~ao 1010, Cidade Universit\'aria, S\~ao Paulo-SP
\newline
Permanent: Saint-Petersburg State University\\
Universitetskaya nab. 7-9, St. Peterburg, Russia\\
e-mail:  wolf86\_666@list.ru}
\end{document}